\documentclass{amsart}
    \usepackage{amsmath}
    \usepackage{amssymb}
    \usepackage[all]{xy}
    \usepackage{comment}
    \usepackage{color}

    \theoremstyle{plain}
    \newtheorem{thm}{Theorem}[section]
    \newtheorem{thm*}{Theorem}[section]
    \newtheorem{cor}[thm]{Corollary}
    \newtheorem{prop}[thm]{Proposition}

    \newtheorem{lemma*}{Lemma}
    
    \newtheorem{recollect}[thm]{Recollection}

    \theoremstyle{definition}
    \newtheorem{defn}[thm]{Definition}
    \newtheorem{remark}[thm]{Remark}

    \newtheorem*{remark*}{Remark}

    \newtheorem{ex}[thm]{Example}
    
    \newtheorem{terminology}[thm]{Terminology}
    
    \newtheorem{question*}{Question}

    \numberwithin{equation}{thm}

    \newcommand{\bB}{\mathbb B}
    \newcommand{\cM}{\mathcal M}
    \newcommand{\cN}{\mathcal N}

    \def\Spec{\operatorname{Spec}\nolimits}

    \def\Proj{\operatorname{Proj}\nolimits}

    \def\stmod{\operatorname{stmod}\nolimits}

    \def\ind{\operatorname{ind}\nolimits}

    \newcommand{\bG}{\mathbb G}

    \newcommand{\bA}{\mathbb A}
    \newcommand{\cA}{\mathcal A}
    
    \newcommand{\cF}{\mathcal F}
    \newcommand{\bP}{\mathbb P}
    \newcommand{\bZ}{\mathbb Z}
    \newcommand{\bF}{\mathbb F}
    
    \newcommand{\cC}{\mathfrak C}
    \newcommand{\cD}{\mathfrak D}
    \newcommand{\cS}{\mathcal S}
    \newcommand{\cK}{\mathcal K}
    
    \newcommand{\cI}{\mathcal I}
    \newcommand{\cE}{\mathcal E}

    \newcommand{\fg}{\mathfrak g}

    \newcommand{\fu}{\mathfrak u}

    \newcommand{\ul}{\underline}

    \def\Spec{\operatorname{Spec}\nolimits}
    \def\sl2{\operatorname{SL_{2(2)}}\nolimits}
    \def\Ga2{\operatorname{\mathbb G_{\rm a(2)}}\nolimits}

    \newcommand{\Gr}{\mathbb G_{(r)}}

    \newcommand{\bN}{\mathbb N}

    \newcommand{\E}{\mathcal E}

    \newcommand{\bU}{\mathbb U}

    \newcommand{\bu}{\bullet}

    \date\today
    
\begin{document}
    
     \title[Support Varieties and Stable Categories]{Support Varieties and Stable Categories for Algebraic Groups}
     
     \author[ Eric M. Friedlander]
    {Eric M. Friedlander$^{*}$\\
     \vskip .2in
     Dedicated to the memory of Brian Parshall} 
    
    \address {Department of Mathematics, University of Southern California,
    Los Angeles, CA}
    \email{ericmf@usc.edu}
    
    \thanks{$^{*}$ partially supported by the Simons Foundation.}

    

    \begin{abstract}  
    We consider rational representations of a connected linear algebraic group $\bG$ over a field
    $k$ of positive characteristic $p > 0$.  We introduce a natural extension \ $M \mapsto \Pi(\bG)_M$ \
    to $\bG$-modules of the $\pi$-point support theory for modules $M$  for a finite 
 group scheme $G$ and show that this theory is essentially 
   equivalent to the more ``intrinsic" and ``explicit" theory \ $M \mapsto \bP\cC(\bG)_M$ \ of supports for an
    algebraic group of exponential type, a theory which uses 1-parameter subgroups 
    $\bG_a  \to  \bG$.
        
    We extend our support theory to bounded complexes of $\bG$-modules, $C^\bu \ \mapsto \ \Pi(\bG)_{C^\bu}$.
   We introduce  the tensor triangulated category $StMod(\bG)$, the Verdier quotient of the bounded
   derived category $D^b(Mod(\bG))$ by the thick subcategory of mock injective modules.  Our support theory 
   satisfies all the ``standard properties"  for a theory of supports for $StMod(\bG)$.
      
   As an application, we employ $C^\bu \ \mapsto \Pi(\bG)_{C^\bu}$ to 
   establish the classification of $(r)$-complete, thick tensor ideals of $stmod(\bG)$
   in terms of  $stmod(\bG)$-realizable subsets of $\Pi(\bG)$ and the classification of  $(r)$-complete, 
   localizing subcategories of $StMod(\bG)$ in terms of 
   $StMod(\bG)$-realizable subsets of $\Pi(\bG)$.
    \end{abstract}
    
    \maketitle
    
   \tableofcontents 
    
    \section{Introduction}
    
    The goal of this work, refining and extending our earlier paper \cite{F2},
    is to present a context and a point of view for the study of
    representations of familiar linear algebraic groups $\bG$ on vector spaces $V$ 
    over a field $k$.  It constitutes a refinement and extension of our earlier work.  We work
    in the modular setting, fixing a field $k$ of prime characteristic $p > 0$; 
    our linear algebraic groups are connected affine group schemes of finite type 
    over $k$; the vector spaces are vector spaces over $k$, not necessarily finite
    dimensional.  A typical example for $\bG$ is the simple algebraic group $SL_n$.
    The representations we consider are ``rational", formally defined as comodules
    for the coordinate algebra $k[\bG]$ with coproduct given by the group 
    structure of $\bG$.
    
    We demonstrate how the theory of support varieties, first developed for finite
    groups and eventually extended to all finite group schemes,  can be extended
    to a good theory of supports for $\bG$-modules.    Our first formulation is an
    evident extension of the theory of $\pi$-points, $M \mapsto \Pi(G)_M$, for 
    $\bG$-modules $M$ with $G$ a finite group scheme as introduced by J. Pevtsova 
    and the author in \cite{FP2}.  This theory,\  \ $M \mapsto \Pi(\bG)_M$, applies to
    any linear algebraic group $\bG$ over an arbitrary field of $k$ of positive 
    characteristic and any (rational) $\bG$-module $M$.  However, it is
    difficult to make explicit because a point of $\Pi(G)$
    is an equivalence class of maps to the group algebra of some Frobenius 
    kernel $\Gr \subset \bG$.  
    Our second formulation $M \mapsto \bP \cC(\bG)_M$, 
    valid only if the algebraic group $\bG$ is of exponential type, uses 
    varieties of 1-parameter subgroups $\bG_a \to \bG$ as introduced in our earlier 
    paper \cite{F2}.

      The underlying justification for a theory of supports is that it is 
 is sensitive to extensions of representations.  
     Most support theories are based on cohomology, and
    support theory offers a geometric picture of cohomology following foundational
    work of D. Quillen \cite{Q}.    For example, 
     the spectrum of the cohomology of the $r$-th Frobenius kernel 
  $(SL_n)_{(r)}$ of $SL_n$ is homeomorphic to the 
    variety of $r$-tuples of pairwise commuting, $p$-nilpotent $n\times n$ 
    matrices of trace 0 \cite{SFB2}.   Our theories are much influenced by
  J. Carlson's consideration of rank varieties for elementary abelian $p$-groups \cite{C1}.  

As motivation for considering support theories for linear algebraic groups, 
    we remind the reader of a few consequences of support theory for 
    finite groups and finite group schemes.
An important role of support theories is that they offer a means of 
    classifying categorical structures associated to representation theories which 
    are ``wild" (e.g., \cite{BCR}, \cite{BIKP}). 
 Support theories have led to the identification and study
    of various interesting special classes of representations such
    as ``modules of constant Jordan type" \cite{CFP} and ``mock injective"
    modules \cite{F4}, provided new invariants for $G$-modules \cite{FPS},  
    and enabled the construction of algebraic vector bundles \cite{FP3}, \cite{BP}.   
    Another fruitful theme has been the application of support theories to 
    investigate structural properties of  the nilpotent cone and related
    algebro-geometric objects \cite{NPV}.
    
   There are obstructions to formulating a suitable support theory for $\bG$-modules.  Not only is 
   rational cohomology of $\bG$ inadequate for this purpose,  but also
   the abelian category $Mod(\bG)$ of $\bG$-modules  rarely has 
   non-trivial projective objects (see \cite{Donk}) and its injective objects are almost 
   always infinite dimensional.   The latter obstruction persuades us to consider 
   the bounded derived category $D^b(Mod(\bG))$ of the abelian category $Mod(\bG)$
   and to formulate associated stable categories.   We verify that our constructions
   $C^\bu \mapsto \Pi(\bG)_{C^\bu}$ and $C^\bu \mapsto \bP \cC(\bG)_{C^\bu}$
   are well defined on objects of the stable category $StMod(\bG)$ and satisfy the 
   standard properties expected of a good theory of supports.
   Although $C^\bu \mapsto \Pi(\bG)_{C^\bu}$ satisfies
   numerous good properties, we are far from recognizing which subspaces of
   $\Pi(\bG)$ are realizable as supports for a bounded complex of $\bG$-modules.
   We lack the analogue for $\bG$-modules of the central result of J. Carlson 
   (easily extended from finite groups to arbitrary finite group scheme $G$) 
   that any closed subvariety of $\Pi(G)$ can be 
   realized as the support of some finite dimensional $G$-module \cite{C2}.
    
    	We highlight some of the contents of this paper.   In Section \ref{sec:Pi-points},
we consider four variations of support theories for modules for a finite group scheme:

 $\quad M \mapsto \bP |G|_M; \quad M  \mapsto \bP V_r(G)_M; \quad  M \mapsto \Pi(G)_M; 
\quad M \mapsto \bP \cC_r(\Gr)_M.$

\noindent
 After briefly recalling the formulations of each, we 
verify in Theorems \ref{thm:p-iso}, \ref{thm:ProjH=Pi}, \ref{thm:1par,pi} and \ref{thm:rho-lambda}
that these theories are equivalent for those finite group schemes $G$ and $G$-modules $M$
for which they are defined.
The $\pi$-point theory $M \mapsto \Pi(G)_M$ applies to any finite group scheme; in Definition
\ref{defn:Pi(bG)}, we show that this theory naturally extends to a theory for linear algebraic groups 
$\bG$ and their $\bG$-modules, $M \mapsto \Pi(\bG)_M$.   

In Section \ref{sec:1-parameter}, we verify for a Frobenius kernel $\Gr$ of a linear algebraic group 
$\bG$ of exponential type that the exponential support theory
$M \mapsto \bP \cC_r(\Gr)_M$ also extends to a theory for $\bG$-modules 
$M \mapsto \bP \cC(\bG)_M$.  In Theorem \ref{thm:cC-Pi}, we show that this
exponential theory $M \mapsto \bP \cC(\bG)_M$ is equivalent to the $\pi$-point theory 
$M \mapsto \Pi(\bG)_M$ for linear algebraic groups of exponential type and arbitrary $\bG$-modules.  
Proposition \ref{prop:further-cC} establishes further good properties of the theory $M \mapsto \bP \cC(\bG)_M$.
In Examples \ref{ex:Ga}, \ref{ex:St-s}. \ref{ex:U3}, and \ref{ex:reductive}
of Section \ref{sec:Gexamples}, we  investigate examples
of $Mod(\bG)$-realizable subsets   $\bP \cC(\bG)_M \ \subset \ \bP \cC(\bG)$.

    For any finite group scheme $G$ over $k$, Theorem \ref{thm:St-rickard}
    extends to arbitrary $G$-modules an important theorem of J. Rickard \cite[Thm2.1]{R1}
    by giving an explicit equivalence of the stable module category $StMod(G)$    
    with the localization of the homotopy category $\cK^b(Mod(G))$ of bounded  complexes of $G$
    by the thick subcategory of bounded complexes of injective $G$-modules.  
    In Definition \ref{defn:stmod}, we introduce the tensor
    triangulated category $StMod(\bG)$ which serves as a natural domain for our support theory 
    for bounded complexes of arbitrary $\bG$-modules.     Motivated by Rickard's equivalence, 
    we formulate our stable module category $StMod(\bG)$ as a localization of the bounded derived category 
   $D^b(Mod(\bG))$ by the thick subcategory of mock injective modules.\
    Defined similarly, the stable module category $stmod(\bG)$ of Definition \ref{defn:stmod} is 
    almost all examples equal to $D^b(mod(\bG))$, the bounded derived category of $mod(\bG)$.
    
    We extend the support theories $M \mapsto \Pi(\bG)_M$ and $M \mapsto \bP \cC(\bG)_M$
    to bounded complexes of $\Gr$-modules and then to bounded complexes of $\bG$-modules
    in Section \ref{sec:Gcomplexes}.
    As seen in Proposition \ref{prop:equalST},  our construction $C^\bu \ \mapsto \Pi(\Gr)_{C^\bu}$ 
     ``agrees" with the usual $\pi$-support theory $M \mapsto \Pi(\Gr)_M$
    for $\Gr$-modules.  
     Theorem \ref{thm:stable-prop} shows that   \ $C^\bu  \ \mapsto  \ \Pi(\bG)_{C^\bu}$
    satisfies all of the standard properties of a theory of supports for $StMod(\bG)$.      Among these properties,
    we mention the criterion for $\Pi(\bG)_{C^\bu}$ to be empty (if and only if $C^\bu$ is an object of $Mock^b(\bG)$),
    the more general fact that $\Pi(\bG)_{C^\bu}$ depends only upon the isomorphism class of $C^\bu$
    in $StMod(\bG)$, and 
   the tensor product property $\Pi(\bG)_{C^\bu\otimes D^\bu} \  = \ \Pi(\bG)_{C^\bu} \cap \Pi(\bG)_{D^\bu}$.

    In Section \ref{sec:stmod(bG)}, we classify certain full subcategories of $stmod(\bG)$ and
    $StMod(\bG)$ in terms of certain subsets of $\Pi(\bG)$.  Namely, Theorem \ref{thm:HPS} establishes
    that $(r)$-complete, thick tensor ideals $\cC \ \subset \ stmod(\bG)$ are classified by locally 
    $stmod(\bG)$-realizable subsets 
    of $\bP (\bG)$, and Theorem \ref{thm:HPS-Mod} demonstrates that $(r)$-complete, localizing  subcategories
    $\tilde \cC \ \subset \ StMod(\bG)$ are classified by $StMod(\bG)$-realizable subsets of $\Pi (\bG)$.  
    The technique of proof of both these theorems
    is to pass from $\bG$ to the family of \ $\{ \Gr, r > 0 \}$ \ of Frobenius kernels of $\bG$.
    For $stmod(\bG)$, we use the classification of thick tensor ideals of $stmod(\Gr)$ given by J. Pevtsova 
    and the author \cite[Thm6.3]{FP2} which in turn utilizes J. Rickard's  idempotents \cite{R2}.  For $StMod(\bG)$, 
    we use the classification of localizing subcategories of $StMod(\Gr)$ proved by 
    D. Benson, H. Krause, S. Iyengar, and J. Pevtsova \cite[Thm10.1]{BIKP}.  
      
    In Section \ref{sec:questions-challenges}, we mention various questions and challenges 
    concerning our support theory which would enhance our understanding of $\bG$-modules.

    Throughout this paper, $k$ will denote a field of characteristic $p > 0$ and all linear algebraic groups
$\bG$ will be assumed to be \ul {connected}; in other words $\bG$ will be assumed to be a 
reduced, connected affine group scheme of finite type over $k$. 
    
    	I am especially indebted to Julia Pevtsova whose suggestion of consideration of colimits
of supports of $\Gr$-modules I finally appreciated.  I also thank Cris Negron for helpful discussions
on matters categorical and Paul Sobaje for patiently answering my questions.   Finally, I 
dedicate this paper to the memory of my friend and collaborator Brian Parshall, who 
contributed much to this mathematics.
	
    \vskip .2in


    \section{Comparison of support theories for finite group schemes}
\label{sec:Pi-points}
    
 We review a few salient features of four formulations of support theories for 
 modules $M$ for a finite group scheme $G$
$$M \mapsto \bP |G|_M; \quad M  \mapsto \bP V_r(G)_M; \quad  M \mapsto \Pi(G)_M; \quad 
M \mapsto \bP \cC_r(\Gr)_M$$ and establish their close relationships.
In Theorem \ref{thm:p-iso}, we recall that the cohomological support theory $M \mapsto \bP |G|_M$ 
agrees with the infinitesimal 1-parameter support theory $M \mapsto \bP V_r(G)_M$ for an
infinitesimal group scheme of height $\leq r$ provided that $M$ is finite dimensional.  In Theorem \ref{thm:ProjH=Pi},
we recall that the cohomological support theory $M \mapsto \bP |G|_M$ is equivalent to the 
$\pi$-point support theory $M \mapsto \Pi(G)_M$ for any finite group scheme provided $M$ is 
finite dimensional.   In Theorem \ref{thm:1par,pi}, we make explicit the equivalence
between $M \mapsto \Pi(\Gr)_M$ and  $M  \mapsto \bP V_r(\Gr)_M$.  

This section concludes with the formulation the exponential support theory $M \mapsto \bP \cC_r(\Gr)_M$
for a Frobenius kernel of a linear algebraic group $\bG$ of exponential type.  Our definition of a linear 
algebraic group of exponential type in Definition \ref{defn:exp-type} allows the exponential 
\ $\cE: \cN_p(\fg) \times \bG_a \ \to \ \bG$ to be a continuous algebraic map (a rational map with unique
specialization at every geometric point).
In Theorem \ref{thm:rho-lambda}, we demonstrate that the exponential support theory 
$M \mapsto \bP \cC_r(\Gr)_M$ agrees with the infinitesimal 1-parameter support theory $M \mapsto \bP V_r(\Gr)_M$
for Frobenius kernels $\Gr$ of linear algebraic groups of exponential type.
This exponential support theory is an an extension of joint work with A. Suslin 
and C. Bendel for infinitesimal group schemes
    \cite{SFB1}, \cite{SFB2} which itself extends earlier joint work with B. Parshall  \cite{FPar}.
    In some sense, $M \mapsto \bP \cC(\bG)_M$ is the ``rank variety" formulation of
    $M \mapsto \Pi(\bG)_M$.

\begin{recollect}
\label{recollect:|G|}
Let $\bG$ be a finite group scheme.
We denote by  \ $H^\bu(G,k)$ \ the commutative, graded $k$-algebra equal to  $H^*(G,k)$ 
if $p = 2$ and otherwise equal to the subalgebra of $H^*(G,k)$ generated by classes of even degree.
We set $\bP |G|$ equal to $\Proj H^\bu(G,k)$.  For $M$ finite dimensional, we denote
by $\bP |G|_M \ \subset \ \bP |G|$ the closed, reduced subscheme determined by 
(the radical of) the homogenous ideal $ker\{ H^\bu(G,k) \to Ext_G^*(M,M) \}$.  
\end{recollect}

The additive group $\bG_a$ plays a central role in what follows.  The coordinate 
    algebra of $\bG_a$ is the polynomial algebra $k[T]$ on one variable, whereas
    the group algebra $k\bG_a$ is the truncated polynomial
algebra $k[u_0,\ldots, u_n, \ldots ]/(u_i^p)$ on countably infinite generators each of whose
$p$-th powers is 0.  Here, $u_r: k[T] \ \to \ k$ is the functional sending
    a polynomial $p(T)$ to its coefficient of $T^{p^r}$.  If we denote by $v_i$
    the dual basis element to $T^i \in k[\bG_a]$ so that $v_i(p(T))$ reads off the coefficient 
    of $T^i$ in $p(T)$, then $u_r = v_{p^r}$.  We denote by
    $\epsilon_r$ the map of $k$-algebras (but not Hopf algebras, unless $r = 1$)
    \begin{equation}
    \label{eqn:epsilon-r}
    \epsilon_r: k\bG_{a(1)} \simeq k[t]/t^p \to k\bG_a, \quad  t \mapsto u_{r-1}.
    \end{equation}
    We also use $\epsilon_r$ to denote the factorization of (\ref{eqn:epsilon-r}) through
    $k\bG_{a(r)} \hookrightarrow k\bG_a$, \ $\epsilon_r: k\bG_{a(1)} \  \to \ k\bG_{a(r)}$.

\vskip .1in

The analysis of $M \mapsto \bP|G|_M$ was achieved for the very special case in which
$G \simeq \bZ/p^{\times r}$ by replacing cohomological support varieties for modules for these
elementary abelian $p$-group by more accessible rank varieties.
This was first generalized to restricted Lie algebras (see \cite{FPar}) and then greatly
expanded to all infinitesimal group schemes in joint work with A. Suslin and C. Bendel
in \cite{SFB1}, \cite{SFB2} as we next recall.

\begin{recollect}
\label{recollect:VrG}
Let $G$ be an infinitesimal group scheme of height $\leq r$.
The functor on commutative $k$-algebras sending $A$ to the set of maps of 
group schemes $\phi: \bG_{a(r),A} \to G_A$
is representable by a graded affine $k$-scheme $V_r(G)$ whose projectivization
we denote by $\bP V_r(G)$.  
For each finite dimensional $G$-modules $M$, one considers the closed subscheme
$\bP V_r(G)_M \ \subset \ \bP V_r(G)$ whose $K$-points for any field extension $K/k$ 
are represented by (infinitesimal)
1-parameter subgroups $\phi: \bG_{a(r),K} \to G_K$ such that $\phi_*(u_{r-1})$ acts 
freely on $M_K$ (in the sense that $(\phi \circ \epsilon_r)*(M_K)$ is a free $K\bG_{a(1)}
\simeq K\bZ/p$-module).
\end{recollect}

\begin{thm}
\label{thm:p-iso}
Let $G$ be an infinitesimal group scheme of height $\leq r$.
There is a natural (in $G$) map
$\psi: H^\bu(G,k) \to V_r(G)$ inducing a universal homeomorphism 
\begin{equation}
\label{eqn:Psi-iso}
\Psi: \bP V_r(G) \ \to \ \bP |G|, \quad \text{restricting to} \quad \Psi: \bP V_r(G)_M \ \to \ \bP |G|_M
\end{equation}
(where $\bP V_r(G)$ denotes $\Proj k[V_r(G)]$) for each finite dimensional $G$-modules $M$.
\end{thm}

The most general construction of a support theory for finite group schemes is the $\pi$-point 
theory \ $M \mapsto \Pi(G)_M$ \
introduced by J. Pevtsova and the author \cite{FP2} which we next recall.

\begin{recollect}
\label{recollect:Pi(G)}
Let $G$ be a finite group scheme $G$ over $k$.  Elements of \ $\Pi(G)$ \
are equivalence classes of flat maps of $K$-algebras
 for some field extension $K/k$,  \ $\alpha_K: K[t]/t^p \ \to \ KG$, \
  which factor through some unipotent abelian subgroup scheme $U_K \subset G_K$.
 Two such flat maps $\alpha_K:K[t]/t^p \ \to \ KG, \ \beta_L:L[t]/t^p \ \to \ LG$ are equivalent if there
 exists a common field extension $\Omega$ of $K, \ L$ such that $\alpha_\Omega^*(M_\Omega)$ is free
 if and only if $\beta_\Omega^*(M_\Omega)$ is free for any finite dimensional $G$-module $M$.
 
 We consider subsets of $\Pi(G)$ of the form
 $$\Pi(G)_M \ \equiv \ \{ [\alpha_K]: \alpha_K^*(M_K)  \ is \ not \ free\ as \ a \ K[t]/t^p{\text-}module \}$$
 for an arbitrary $\bG$-module $M$, well defined by \cite[Thm6.6]{F2}.  The closed subsets of $\Pi(G)$ are the 
 subsets of this form with $M$ finite dimensional.
\end{recollect}

The following theorem tell us us that $M \mapsto \bP |G|_M$ is equivalent to $M \mapsto \Pi(G)_M$
provided that $M$ is finite dimensional.

\begin{thm} \cite[Thm 7.4]{FP2}
\label{thm:ProjH=Pi}
There exists a natural scheme structure on $\bP |G|$ and a scheme-theoretic isomorphism
\begin{equation}
\label{eqn:Phi-iso}
\Phi:  \bP |G| \ \stackrel{\sim}{\to} \ \Pi(G), \quad \text{restricting to} \quad \bP |G|_M 
\ \stackrel{\sim}{\to} \ \Pi(G)_M
\end{equation}
for any finite dimensional $G$-module $M$.
\end{thm}

The isomorphisms $\Psi$ of (\ref{eqn:Psi-iso}) and  $\Phi$ of (\ref{eqn:Phi-iso}) are 
quite abstract.  The natural map $\psi: H^\bu(G,k) \ \to \ k[V(G)]$ of \cite[Thm1.14]{SFB1} 
determined $\Psi$ for an infinitesimal group scheme of height $\leq r$ uses the universal height $r$
 1-parameter subgroup for $G$.   The isomorphism $\Phi$ of \cite[Thm7.5]{FP2} entails 
 consideration of a sheafified version of the endomorphisms of the stable module 
 category of finite dimensional $G$-modules.
 
 The following proposition makes $\Psi, \ \Phi$ more concrete.

\begin{thm}
    \label{thm:1par,pi}
Let $G$ be an infinitesimal group scheme over $k$ of height $\leq r$.
Then the composition
$$\bP V_r(G) \quad \stackrel{\Psi}{\to} \quad \bP |G| \quad \stackrel{\Phi}{\to} \Pi(G)$$
sends the $K$-point represented by the 1-parameter subgroup $\eta_K: \bG_{a(r),K} \to G_K$ to the point
 represented by the $\pi$-point
\ $\eta_{K*} \circ \epsilon_r: K[t]/t^p \to K\bG_{a(r),K} \to KG_K$.
\end{thm}

 \begin{proof}
 By \cite[Thm3.6]{FP2}, the inverse of $\Phi$ (as a homeomorphism) sends a $\Pi$-point 
 $\alpha_K: K[t]/t^p$ to 
 \begin{equation}
 \label{eqn:ker1}
 ker\{ \alpha_K^*: H^\bu(G_K,K) \to H^\bu(\bG_{a,K},K\} \cap H^\bu(G,k).
 \end{equation}
 
 To identify $\Psi$, we first consider the special case $G = GL_{N(r)}$.  
 By \cite[Prop2.9]{FP1}, each equivalence class of $\pi$-points of $GL_{N(r)}$ 
 contains a representative of the form 
 $$ \alpha_K \ = \ \cE_{\ul B*} \circ \epsilon_r: K[T]/t^p \to K\bG_{a(r),K}  \to k\bG_K,$$ 
 unique up to base change by $L/K$ and up to scalar multiples.   Then $\Psi$ sends a
 $K$-point $k[V_r(G)] \to K$ given by evaluation at the 1-parameter subgroup $\cE_{\ul B}$
 defined over $K$ to 
 \begin{equation}
 \label{eqn:ker2}
 ker\{ ev_{\ul B} \circ \psi: H^\bu(G,k) \to k[V_r(G)] \to K\}.
 \end{equation} 
 Thus, for $G = GL_{n(r)}$, it suffices to show the equality of the radicals of the kernels of 
 (\ref{eqn:ker2}) and of 
 \begin{equation}
 \label{eqn:ker3}
 ev_{\ul B} \circ \psi: H^\bu(G,k) \to k[V_r(G)] \to K, \quad
 (\cE_{\ul B*} \circ \epsilon_r)^*: H^\bu(G_K,K) \to H^\bu(\bG_{a,K},K\} \cap H^\bu(G,k).
 \end{equation} 
 This follows from \cite[Thm5.2]{SFB1}.

 For a general infinitesimal group scheme $G$, we embed $G$ in some $GL_N$ and use \cite[Thm5.2]{SFB2}
 which asserts that the map $\psi: H^\bu(G,k) \to k[V_r(G)]$ has nilpotent kernel and
 \cite[Prop4.3]{SFB2} which implies that $H^\bu(GL_N,k) \to H^\bu(G,k)$ has image containing all $p^r$-th powers
 so that \ $\bP |G|  \ \to \ \bP |GL_{N(r)}|$ \ is injective.
 We consider the diagram
  \begin{equation}
 \label{eqn:GLn-G}
 \begin{xy}*!C\xybox{%
   \xymatrix{ 
\bP V_r(GL_{N(r)}) \ar[r]^\Psi \ar[d] & \bP |GL_{N(r)}| \ar[r]^\Phi \ar[d] & \Pi(GL_{N(r)}) \ar[d] \\
 \bP V_r(G) \ar[r]^\Psi & \bP |G| \ar[r]^\Phi  & \Pi(G)  }
    }\end{xy}
 \end{equation}
 which commutes by the naturality of $\psi$ (which determines $\Psi$) and $\Phi$.
 Thus, a diagram chase 
 verifies that the radicals of the two kernels of the maps of algebras (\ref{eqn:ker1}), (\ref{eqn:ker2}) are equal.
  \end{proof}

\vskip .1in

The abstract formulation of $M \mapsto \Pi(G)_M$ does not easily lead to computations.  
Following \cite{SFB1}, we consider the following affine $k$-varieties, worthy of study
in their own right.

   \begin{defn}
    \label{defn:NpG}
    Let $\bG$ be a linear algebraic group with Lie algebra $\fg$.  We
    denote by $\cN_p(\fg) \subset \fg$ the $p$-nilpotent variety of $\fg$, the 
    (reduced) closed subvariety of the affine space $\Spec S(\fg^*) \simeq \bA^d$ 
    whose $K$-points are the elements $X \in \fg_K$ such that $X^{[p]} = 0$ for any field  
    extension $K/k$.     For any $r > 0$, we define 
    $$\cC_r(\bG) \ \equiv \ \cC_r(\cN_p(\fg)) \ \subset \ \cN_p(\fg)^{\times r}$$
    to be the (reduced) closed subvariety of $\cN_p(\fg)^{\times r}$ whose $K$-points
    are $r$-tuples $\ul B = (B_0,\ldots,B_{r-1}) \in \cN_p(\fg_K)^{\times r}$ satisfying the
    condition the $[B_i,B_j] = 0$ for all $i,j$ with $0 \leq i < j < r$. 
 \end{defn}

      The following definition of a linear algebraic group of exponential type
    is a slight modification of that of \cite{F2} in that we allow $\cE$ to be a continuous algebraic 
    map and we explicitly require $\cE$  to be $\bG$-equivariant.   Following \cite{F1}, we define a 
    continuous algebraic map
$f: X \ \dashrightarrow \ Y$ to be a ``roof" $X \stackrel{p}{\leftarrow} \tilde X \stackrel{g}{\to} Y$
of $k$-schemes such that  $p: \tilde X \to X$ is a universal homeomorphism (i.e., a finite, surjective map
such that $k(\tilde x)$ is purely inseparable of $k(p(\tilde x))$ for all points $\tilde x \in \tilde X$).  A typical example
    is a rational map from $X$ to $Y$ (i.e., a morphism with domain a dense open set) whose graph
    in $X\times Y$ projects to $X$ via a universal homeomorphism.   A bicontinuous algebraic map of irreducible
 varieties $X, \ Y$, 
$f: X \stackrel{\sim}{\dashrightarrow} Y$, is given by a finite, purely separable extension $k(X)$ over $k(Y)$
inducing a bijection between the $K$-points of $X$ and the $K$-points of $Y$ for any algebraically closed
extension $K$ of $k$.  Examples include the Frobenius map $F: \bG \to \bG^{(1)}$ and the normalization
$\tilde N(sl_2) \to N(sl_2)$.
    
    	As remarked in \cite[Rem1.7]{F2}, if the linear algebraic group $\bG$ admits the structure 
$(\bG,\cE)$ of an algebraic group of exponential type, two such structures are isomorphic via 
an automorphism of $\cN_p(\fg)$.
           
    \begin{defn}
    \label{defn:exp-type}
    Let $\bG$ be a linear algebraic group over $k$ with Lie algebra $\fg$ equipped with a 
    $\bG$-equivariant, continuous algebraic map $\cE: \cN_p(\fg) \times \bG_a \ \dashrightarrow  \ \bG$ 
    sending a geometric point $(B,\alpha)$ of $\cN_p(\fg) \times \bG_a$ to $\cE_B(\alpha)$.  
    Then $(\bG, \cE)$ is said to be an {\it algebraic group of exponential type} 
    provided that 
    \begin{enumerate}
    \item
    For each Frobenius kernel $\Gr \ \subset \bG$ and geometric point $\Spec K \to V(\Gr)$,
    the corresponding 1-parameter subgroup  $\bG_{a(r),K} \ \to \ \bG_{(r),K}$ admits a unique representation 
    of the form 
    $$ \cE_{\ul B}  \ \equiv \ \large\prod_{s = 0}^{r-1} \cE_{B_s} \circ F^s: \bG_{a(r),K} \ \to \ \bG_{(r),K}.$$
    for some $K$-point 
    $\ul B = (B_0,\ldots,B_{r-1}) \in \cC_r(\bG).$
    Here, $F: \bG_a \to \bG_a$ is the Frobenius map, and 
    $\cE_B: \bG_{a,K} \to \bG_K$ is the 1-parameter subgroup associated to the map of $K$-points
    determined by $\cE$  restricted to $\{ B \} \times \bG_a$.
\item
    For each algebraically closed field extension $K/k$, a 1-parameter 
    subgroup $\bG_{a,K} \to \bG_K$   has a unique representation of the form
    $$ \cE_{\ul B}  \ \equiv \ \large\prod_{s\geq 0}\cE_{B_s} \circ F^s:  \bG_{a,K} \ \to \ \bG_{K}.$$
    for some $K$-point $\ul B = (B_0,\ldots,B_n,\ldots)$ of 
    $\cC(\bG) \ \equiv \ \varinjlim_r \cC_r(\bG)$.
 \end{enumerate}
  \end{defn} 
  
  \vskip .1in
  
  We are much indebted to P. Sobaje for guiding us to the following results in the literature.

\begin{thm}
\label{thm:exp-type}
The following are examples of linear algebraic groups $\bG$ of exponential type:
\begin{enumerate}
\item
Simple algebraic groups $\bG$ of classical type, their standard parabolic subgroups, 
and the unipotent radicals of these standard parabolic subgroups.   See \cite{SFB1}.
\item 
Connected, reductive groups $\bG$ over an algebraically closed field $k$ whose derived subgroup has no factor of 
exceptional type with Coxeter number $h$ with $p \leq  2h-2$ .  See \cite[Thm9.5]{McN} for a somewhat 
sharper result, with an explicit list of problematic exceptional simple factors.
\item
$\bU$\ is the unipotent radical of a parabolic subgroup $\bP \subset \bG $ of a connected, reductive group  
over an algebraically closed field $k$ with the property that the nilpotent class of $\bU$ is $< p$.
See  \cite[Prop5.3]{Sei}.
\end{enumerate}
\end{thm}

\begin{proof}
The examples of (1) are given in \cite[Lem1.8]{SFB1}, where $\cE: \cN_p(\fg) \times \bG_a \to \bG$ 
is constructed 
as a scheme-theoretic morphism and no condition is placed on the field $k$.

In the context of (2), G. McNinch establishes $\cE: \cN_p(\fg) \times \bG_a \to \bG$ as an ``isomorphism of varieties"
(but not of schemes).  Moreover, $\cE$ is given by the restriction of the exponential $exp: \cN_p(gl_N) \times \bG_a \to GL_N$
associated to some faithful representation $(\rho,V)$ of $\bG$ of some dimension $N$.  Although McNinch only considers infinitesimal
1-parameter subgroups, his results also apply to establish condition (2) of Definition \ref{defn:exp-type} 
when supplemented by the following observation.  Given a 1-parameter subgroup $\phi: \bG_a \to \bG$, 
condition (1) of Definition \ref{defn:exp-type} tells us that $(\rho \circ \phi)_{|\Gr}$ is given uniquely by a 
product of exponentials; using 
the observation of \cite{SFB1} that 1-parameter subgroups $\bG_a \to GL_N$ have a unique, finite product representation,
we conclude that the product representations for $(\rho \circ \phi)_{|\Gr}$ agree for sufficiently large $r$ and thus
that $\rho \circ \phi$ agrees with this ``stable" product representation of $(\rho \circ \phi)_{|\Gr}$.

In the context of (3), G. Seitz proves in \cite[Prop5.3]{Sei} that there is a unique $\bP$-equivariant isomorphism
$\theta: \fu \ \stackrel{\sim}{\to} \ \bU$ of algebraic groups, where $\fu$ is viewed as a vector group over $k$.  Thus, any
1-parameter subgroup $\bG_a \to \bU$ is given by composing an additive map $\bG_a \to \fu$ with $\theta$.
\end{proof}

In the following definition, we consider $p$-nilpotent linear operators on $K[\bG]$ of the form \ $\cE_{B*}(u_r)$.
We recall that $u_r:k[\bG_a] \to k$ reads of the coefficient of $T^{p^{r-1}}$ of a polynomial in in $k[T] = k[\bG_a]$
and $\cE_{B*}(u_r)$ is the composition \ $u_r \circ \cE_B^*: K[\bG] \to K[\bG_a] \to K$.

  \begin{defn}
\label{defn:geom-Cr}
Let $(\bG,\cE)$ be a linear algebraic group of exponential type.  For any  $\bG$-module $M$, 
we define the {\it Jordan type} of $M$ at a $K$-point $\cE_{\ul B}$ of $\cC_r(\cN_p(\fg))$
to be block sum decomposition 
 $$JT(M)_{\cE_{\ul B}} \ \equiv  \ [1]^{\oplus m_1} \oplus \cdots \oplus [p]^{m_p}, \ 0 \leq m_i \leq \infty$$
for the action on $M_K$ of the $p$-nilpotent element $\sum_{s \geq 0} (\cE_{B_s})_*(u_s) \in K\bG;$
we further define the ``stable Jordan  type of $M$ at $\cE_{\ul B}$ to be 
$$sJT(M)_{\cE_{\ul B}} \ \equiv  \ [1]^{\oplus m_1} \oplus \cdots \oplus [p-1]^{m_{p-1}}, \ 0 \leq m_i \leq \infty.$$

We define $\cC_r(\cN_p(\fg))_M \ \subset \ \cC_r(\cN_p(\fg))$ to be the subspace whose set of $K$-points 
is the subset of those $K$-points $\ul B$ of $\cC_r(\cN_p(\fg))$ with the property that the stable Jordan type
$sJT(M)_{\cE_{\ul B}}$ is not 0.  

Equivalently (see\cite[Defn4.4]{F2}), the $K$-point $\ul B$ is a $K$-point of $\cC_r(\cN_p(\fg))_M$
\ $\iff$ \ $(\cE_{\Lambda_r(\ul B)}) \circ \epsilon_r)^*(M_K)$ is not free as a $K[t]/t^p$-module.
\end{defn}

\vskip .1in

 The natural grading on the affine scheme $V_r(\Gr)$ for $(\bG,\cE)$ a linear
 algebraic group of exponential type corresponds to the monoid action
 $$V_r(\Gr) \times \bA^1 \ \to \ V_r(\Gr), \quad (\cE_{\ul B}, a) \ \mapsto  \cE_{a \bu \ul B},$$
 where $\cE_{a \bu \ul B} \ = \  \E_{\ul B}(a \cdot t)$.
  Since $(\cE_B\circ F^s)(a\cdot t) = \cE_{a^s \cdot B}(t)$,
 we conclude that $V(\Gr) \times \bA^1 \ \to \ V(\Gr)$ is given by
 \begin{equation}
 \label{eqn:A1-action}
  (\cE_{\ul B},\alpha) \ \mapsto \ \cE_{(\alpha\bu \ul B)}, \quad \alpha \bu (B_0,\ldots B_{r-1}) \equiv
  (a\cdot B_0,\ldots,a^{p^{r-1}}\cdot B_{r-1}).
 \end{equation}
 We denote by $\bP V_r(\Gr)$ the projective variety associated to $V_r(\Gr)$.
 
 We introduce a grading on $\cC_r(\bG)$ which enables 
 bicontinuous algebraic maps  $\rho_r \circ \Lambda_r: \bP \cC_r(\bG) \to \ \bP V_r(\Gr)$.
 Namely, we define the  monoid action
\ $\cC_r(\bG)  \times \bA^1 \ \to \ \cC_r(\bG) $ \
as the restriction of the monoid action \ $\fg^{\times r} \times \bA^1 \to \fg^{\times r}$
given by
 \begin{equation}
 \label{eqn:A*-action}
  (\ul B,a) \ \mapsto \ a * \ul B, \quad  a* (B_0,\ldots B_{r-1}) \equiv
  (a^{p^{r-1}}\cdot B_0,\ldots,a \cdot B_{r-1}).
 \end{equation}
 We denote by $\bP \cC_r(\bG)$ the associated projective variety.  For any $\bG$-module
 $M$ and any algebraically closed field extension $K/k$, we define the $K$-points of 
 $\bP\cC_r(\bG)_M$ to be the image 
 $K$-points of $\cC_r(\bG)_M \backslash \{ 0 \}$. 

 \vskip .1in
 
 \begin{defn}
 \label{defn:rho}
 Let $(\bG,\cE)$ be a linear algebraic group of exponential type.   We denote by 
 \begin{equation}
\label{rho-r}
\rho_r:  \cC_r(\bG) \quad \dashrightarrow \quad V_r(\Gr).
\end{equation}
the continuous algebraic map associating to 
a $K$ point  $\ul B$ of $\cC_r(\bG)$  the $K$-point $\cE_{\ul B}: \bG_{a(r),K} \to \bG_{(r),K}$
 of $V_r(\Gr)$; so defined, $\rho_r$ is a 
 bicontinuous algebraic map and thus a universal homeomorphism.
 Here, we have abused notation by using $\cE_{\ul B}$  to denote both the 1-parameter subgroup 
 $\bG_{a,K} \to \bG_K$ and its  restriction $(\cE_{\ul B})_{(r)}: \bG_{a(r),K} \to \bG_{(r),K}$.
 
 We denote by \ $\Lambda_r:  \cC_r(\bG) \to \ \cC_r(\bG)$ \
 the isomorphism sending $(B_0,\ldots,B_{r-1})$ to $(B_{r-1},\ldots,B_{0})$.
\end{defn}

\vskip .1in

 \begin{thm}
 \label{thm:rho-lambda}
Let $(\bG, \cE)$ be an algebraic group of exponential type and let $M$ be a $\bG$-module.  
The bicontinuous algebraic map given as the composition
\begin{equation}
    \label{eqn:rho-lambda}
   \rho_r \circ \Lambda_r:\cC_r(\bG) \ \stackrel{\sim}{\to} \  \cC_r(\bG) \ \stackrel{\sim}{\dashrightarrow} \ V_r(\Gr)
    \end{equation}
    commutes with the gradings of $\cC_r(\bG)$ given by (\ref{eqn:A*-action}) and of $V_r(G)$ given
    by (\ref{eqn:A1-action}).
    
    Furthermore, $\rho_r \circ \Lambda_r$ satisfies the following properties:
\begin{enumerate}
 \item
For all $\bG$-modules $M$, $\rho_r \circ \Lambda_r$ determines a bicontinuous algebraic map     
\begin{equation}
    \label{eqn:rho-lambda-r}
\rho_r \circ  \Lambda_r: \bP\cC_r(\bG)  \ \stackrel{\sim}{\dashrightarrow}   \ \bP V_r(\Gr), 
\quad \text{restricting to} \ \ \bP\cC_r(\bG)_M  \ \stackrel{\sim}{\dashrightarrow}  \ \bP V_r(\Gr)_{M_{|\Gr}}
    \end{equation}
    where $M_{|\Gr}$ denotes the restriction of $M$ to $\Gr$.
\item
If $M$ is finite dimensional, then $\rho_r \circ \Lambda_r$ restricts to a universal homeomorphism from
$ \bP\cC_r(\bG)_M \ \subset \ \bP\cC_r(\bG)$ \ to \ $\bP V_r(\Gr)_{M_{|\Gr}}
 \ \subset \ \bP V_r(\Gr)$.
\item
The adjoint action of $\bG$ on $\fg$
 induces a $\bG(K)$-action on the $K$-points of $\bP \cC(\bG)$ for any field extension $K/k$.
 This action stabilizes the $K$-points of  $\bP\cC_r(\bG)_M$ for any $\bG$-module $M$.
 \item
 The composition $\Psi \circ (\rho_r \circ  \Lambda_r): \bP \cC_r(\bG) \ \simeq \ \bP |\Gr|$
 is $\bG(K)$ equivariant when evaluated at $K$-points for any field extension $K/k$.
 \end{enumerate}
\end{thm}

\begin{proof}
Comparing the actions of $\bA^1$ on $\cC_r(\bG)$ as in (\ref{eqn:A*-action}) and of $\bA^1$ on 
$V_r\bG)$ as in (\ref{eqn:A1-action}), we see that $\rho_r \circ \Lambda_r$ respects degrees of
homogeneous elements.

By \cite[Prop4.3]{F2}, the $\pi$-point  $K[u]/u^p \to K\Gr$ given by $u \mapsto \sum_{s=0}^{r-1} (\cE_{B_s})_*(u_s)$ 
is equivalent to the  $\pi$-point sending $u \mapsto (\cE_{\Lambda_r(\ul B)})_*(u_{r-1})$ provided that $B_s = 0 $ for $s \geq r$.
Thus, the Jordan type of  $M$ at the geometric point $\cE_{\ul B} $ of $\cC_r(\cN_p(\fg))$ 
as defined in Definition \ref{defn:geom-Cr} has a block of size $< p$ consists 
 if and only if the 1-parameter subgroup $\cE_{\Lambda_r(\ul B)}: \bG_{a(r),K} \to \bG_{(r),K}$ lies in $V_r(\Gr)_M$.
This implies that (\ref{eqn:rho-lambda-r}) induces a bijection on $K$ points for any $K/k$ and thus is a 
bicontinuous map.

To prove that $\bP \cC_r(\bG)_M \subset \ \bP\cC_r(\bG)$ is closed if $M$ is finite dimensional, 
we observe that the condition on a  $K$-point $\ul B$ of $\bP \cC_r(\bG)$ to lie
in $\bP \cC_r(\bG)_M$ is the condition that the rank of
$(\cE_{\Lambda_r(\ul B)})_*(u_{r-1}) \in K\bG$ acting on $M_K$
is strictly less that $\frac{p-1}{p} dim(M)$.    

Observe that $(M_K)^\tau$ is isomorphic to $M_K$ as a $\bG$-module for any $\bG$-module
$M$ and any $\tau \in \bG(K)$.   Moreover, for any $\tau \in \bG(K)$, the
 uniqueness argument of \cite[Remark1.7]{F2} and the $\bG$-equivariance of Springer isomorphisms
 imply that  $(\cE_{\ul B})^\tau  \ = \ \cE_{\ul B^\tau}$.    Thus, $((\cE_{\ul B})^\tau)^*(M_K)$ is isomorphic
 to $(\cE_{\ul B})^*(M_K)$, thereby proving (3).
 
Finally, $\Psi \circ (\rho_r \circ  \Lambda_r)$ sends a $K$-point $\ul B$ of $\cC_r(\bG)$ to the 
intersection with $H^\bu(G,k)$ with the kernel of 
$((\cE_{\Lambda_r(\ul B)})_*\circ \epsilon_r)^*: H^\bu(\Gr,K) \to H^\bu(K[t]/t^p,K)$.  For any $\tau \in \bG(K)$,
\ $\tau \circ \cE_{\ul B} \ = \ \cE_{\ul B^\tau}$ \ by $\bG$-equivariance for $\cE$.  Thus, $\tau$ applied to
this kernel equals the intersection with $H^\bu(G,k)$ with the kernel of 
$(\cE_{\Lambda_r(\ul B^\tau)})_*\circ \epsilon_r)^*$

\end{proof} 
\vskip .1in

    \vskip .2in
    

\section{Support theories $M \mapsto \Pi(\bG)_M, \ M \mapsto \cC(\bG)_M$}
    \label{sec:1-parameter}
    
In this section, we consider support theories for $\bG$-modules for a linear algebraic group.
Our first theory, $M \mapsto \Pi(\bG)_M$ is a natural extension of the $\pi$-point support
theory for finite group schemes recalled in Recollection \ref{recollect:Pi(G)}.   Although simple
to define and good for establishing general properties, $M \mapsto \Pi(\bG)_M$ does not lend itself to
computation.  With this in mind, we also consider the natural extension $M \mapsto \bP\cC(\bG)_M$
of the exponential theory for Frobenius kernels of linear algebraic groups of exponential 
type formulated in Definition \ref{defn:NpG}.

For notational convenience, we shall frequently denote the restriction $M_{|\Gr}$ 
 of a $\bG$-module $M$ to some Frobenus kernel$\Gr \ \subset \ \bG$ by $M$.

	We utilize $M \mapsto \bP\cC(\bG)_M$ to verify various properties of $M \mapsto \Pi(\bG)_M$
for $\bG$ a linear algebraic group of exponential type.  We also use  $M \mapsto \bP\cC(\bG)_M$ 
in order to consider $\bG$-modules of bounded 
exponential degree in Proposition \ref{prop:exp-deg}.

	We begin with an observation about closed embeddings of infinitesimal group schemes.
This observation is contrary to the behavior of cohomological support varieties for finite groups.
For example, if $\tau$ is a finite group and $\kappa \subset \tau$ is a $p$-Sylow subgroup, then
$|\kappa| \to |\tau|$ is rarely injective.
    
    \begin{prop}
 \label{prop:pi-functor}
 Let $i: H \hookrightarrow G$ be a closed embedding of infinitesimal group schemes of 
 height $\leq r$.  Then
 for any $G$-module $M$, $i$ induces an embedding \quad $i_*: \Pi(H) \ \subset \ \Pi(G)$ which
 restricts to $\Pi(H)_{i^*M}\subset \ \Pi(G)_M$.
 \end{prop}
 
 \begin{proof}
 Clearly, composition with $i$ induces an embedding $V_r(H) \to V_r(G)$ and thus by
 Theorem 1.6 an embedding $i: \Pi(H) \to \Pi(G)$ also given by composition with $i$.
 In other words, sending a $\pi$-point $\alpha_K: K[t]/t^p \to KH$ to $i_*\circ \alpha_K: 
 K[t]/t^p \to KH \to KG$ is well defined and injective on equivalence classes of $\pi$-points.
 (Recall that $i_*: kH \to kG$ is flat whenever $i: H \to G$ is a closed embedding of
 finite group schemes, so that $i\circ \alpha_K$ is always a 
 $\pi$-point whenever $\alpha_K$ is a $\pi$-point for any closed embedding of finite
 group schemes; see \cite[8.16.2]{J}.) 
 
 Let $M$ be a $G$-module and $\alpha_K$ be a $\pi$-point of $H$.  Then $i\circ \alpha_K$
 represents a point in $\Pi(G)_M$ $\iff$ $(i\circ \alpha_K)^* (M_K)$ is not free
 as a $K[t]/t^p$-module $\iff$ $\alpha_K^*(i^*M_K)$ is not not free as a $K[t]/t^p$-module 
 $\iff$ $\alpha_K$ represents a point of $\Pi(H)_{i^*M}$.  Thus, $i_*: \Pi(H) \to \Pi(G)$ restricts to 
$\Pi(H)_{i^*M} \to \Pi(G)_M$ for any $G$-module $M$. 
 \end{proof}

 Proposition \ref{prop:pi-functor} justifies the constructions of 
 the following definition of $M \mapsto \Pi(\bG)_M$ as a colimit
with respect to $r$ of $M \mapsto \Pi(\Gr)_M$.  One can view $M \mapsto \Pi(\bG)_M$ as a
support theory for modules for the hyperalgebra $k\bG \ \equiv \ \varinjlim_r \Gr$.
    
      \begin{defn}
  \label{defn:Pi(bG)}
  Let $\bG$ be a linear algebraic group over $k$.  We define $\Pi(\bG)$ to be the topological space 
  \quad $\Pi(\bG) \quad \ \equiv \quad \varinjlim_r \Pi(\Gr)$, \ the colimit with the colimit topology
  whose connecting maps $ \Pi(\Gr) \to \Pi(\bG_{(r+1)})$  are given by sending a $\pi$-point
  $\alpha_K: K[t]/t^P \to K\bG_{(r)}$ to its composition with the flat map $K\bG_{(r)} \to K\bG_{(r+1)}$
  induced by the closed embedding $\Gr \hookrightarrow \bG_{(r+1)}$.
  
  For a $\bG$-module $M$,   the $\Pi$-support space $\Pi(\bG)_M$ is defined to be 
  $$\Pi(\bG)_M \ \equiv \ \varinjlim_r \Pi(\Gr)_{M_{|\Gr}}
  \quad \ \subset \quad \varinjlim_r \Pi(\Gr) \ \equiv \Pi(\bG).$$
  \end{defn}
  \vskip .1in
  The fact that $\Pi(\Gr) \hookrightarrow \Pi(\bG_{(r+1)})$ restricts to 
  $\Pi(\Gr)_{M_{|\Gr|}} \hookrightarrow \Pi(\bG_{(r+1)})_{M_{|\bG_{(r+1)}|}}$
for $M$ a finite dimensional $\bG_{(r+1)}$-module is immediate from the
  definition of the equivalence relation on $\pi$-points.  For arbitrary $\bG_{(r+1)}$-modules $M$,  one appeals
  to \cite[Thm4.6]{FP2} which asserts that equivalence of $\pi$-points implies strong equivalence.

 \begin{remark}
 We compare  our current construction of $M \to \Pi(\bG)_M$  with the construction of 
 $M \to V(\bG)_M$ considered in \cite[Defn4.4]{F2}.
First of all, \cite{F2} requires that the linear algebraic group $\bG$ be of exponential type. 
Second, in \cite{F2}, $k$ is required to 
be algebraically closed and the formulation of support varieties is as a topological 
space of closed points, whereas
$\Pi(\bG)_M$ includes non-closed points.
Finally, $\Pi(\bG)_M$ as defined in Definition \ref{defn:Pi(bG)} is a colimit with respect to $r$ 
of $\Pi(\Gr)_{M_{|\Gr|}}$, whereas limits were taken in \cite{F2} rather than colimits. 
\end{remark}
     
We recall the definition of a mock injective module for a linear algebraic group $\bG$ introduced
in \cite[Defn4.3]{F2}..  The
subcategory of $Mod(\bG)$ consisting of mock injective modules plays a key role in our formulation
 of stable module categories in Section \ref{sec:stab-cat}.   Interesting examples of such modules,
    even for $\bG_a$, are constructed in \cite{F2} and \cite{HNS}.

 \begin{defn}
 \label{defn:mock}
 Let $\bG$ be a linear algebraic group.  A $\bG$-module $J$ is said to be mock injective if 
 the restriction of $J$ to each Frobenius kernel $\Gr$ of $G$ is an injective $\Gr$-module.
 \end{defn} 
  
    The properties for our support theory $M \mapsto \Pi(\bG)_M$ 
    stated in the following theorem are the basic properties required of a 
    theory of support.   Perhaps it is worth mentioning that support theories do not determine
    functors from categories of $G$-modules to subsets of some space.  For example, 
    one could consider a $\bG$-module $M$ with non-trivial support and an
    embedding $M \hookrightarrow I$ of $M$ into an injective $\bG$-module; 
    in this case, there is no conceivable map from the support of $M$ to 
    the empty set  (which is the support of $I$).
    
      Since $M \to \Pi(\bG)_M$ is the colimit of $M \to \Pi(\Gr)_{M_{|Gr}}$, 
the properties of $M \mapsto \Pi(\bG)_M$ 
stated in Theorem \ref{thm:properties} are immediate consequences 
of the corresponding properties for each Frobenius kernel $\Gr$ of $\bG$
as established in \cite{FP2} together with the definition of Mock injective modules.

    \begin{thm}
    \label{thm:properties}
    Let $(\bG, \cE)$ be a linear algebraic group and consider $\bG$-modules $M, M_i, N$.
    Then sending a $\bG$-module $M$ to the subspace \ $\Pi(\bG)_M \subset \Pi(\bG)$ \ satisfies
    the following properties:
    \begin{enumerate}
    \item
    {\it Isomorphism}: If $M$ and $N$ are isomorphic $\bG$-modules, then (as subsets of $\Pi(G)$)
    $$\Pi(G)_M \quad = \quad \Pi(G)_N.$$
    .
    \item
    {\it Arbitrary direct sums}: For any family $\{ M_i, i \in I\}$ of $\bG$-modules, 
    $$\Pi(\bG)_{\oplus_I M_i} \quad = \quad \bigcup_{i\in I} \Pi(\bG)_{M_i}.$$
    \item
    {\it Tensor products}: For any pair of $\bG$-modules $M$ and $N$, 
    $$\Pi(\bG)_{M\otimes N} \ = \ \Pi(\bG)_M \cap V(\bG)_N.$$
    \item
    {\it Two out of three}:
    For any short exact sequence of $\bG$-modules $0 \to M_1 \to M_2 \to M_3 \to 0$ and any
    permutation $\sigma$ of $\{ 1, 2, 3 \}$, 
    $$\Pi(\bG)_{M_{\sigma(2)}} \ \subset \ \Pi(\bG)_{M_{\sigma(1)}} \cup \Pi(\bG)_{M_{\sigma(3)}}.$$
    \item
    {\it Trivial module}: \quad $\Pi(\bG)_k \quad = \quad \Pi(\bG)$.
    \item
    {\it Closed}: If $M$ is finite dimensional, then $\Pi(\bG)_M\  \subset \ \Pi(\bG)$ is a closed subspace.
    \item
    {\it Detection}: $\Pi(\bG)_M = \emptyset$ $\quad \iff \quad$ $M$ is mock injective.
    \end{enumerate}
    \end{thm}

   \vskip .1in
   
To establish further properties and to compute examples, we shall consider the colimit of the exponential
  support theory $M \mapsto \bP \cC_r(\bG)_M$ of Definition \ref{defn:geom-Cr}.  
We observe that the natural embedding $\cC_r(\bG) \ \hookrightarrow \ \cC_{r+1}(\bG)$ (determined by
sending a $K$-point $\ul B = (B_0,\ldots,B_{r-1})$ to $(\ul B,0) \equiv (B_0,\ldots,B_{r-1},0)$ multiplies
the degrees of homogeneous elements (specified by the $\bA^1$-action given in (\ref{eqn:A*-action}))
by $p$ and thus induces $\bP \cC_r(\bG) \ \hookrightarrow \ \bP \cC_{r+1}(\bG)$.  Moreover, for 
$(\bG,\cE)$ an algebraic group of exponential type, one checks by
inspection the equality for any $K$-point $\ul B$ of $\cC_r(\bG)$ 
\begin{equation}
\label{eqn:epsilons}
\cE_{\Lambda_r(\ul B)} \circ \epsilon_r \quad = \quad \cE_{\Lambda_{r+1}((\ul B,0))} \circ \epsilon_{r+1}: K[t]/t^p \to \ \bG.
\end{equation}

 \vskip .1in
 \begin{thm}
 \label{thm:cC-Pi}
 Let $(\bG, \cE)$ be an algebraic group of exponential type and consider the 
 bicontinuous algebraic map $ \rho_r \circ \Lambda_r:  \bP\cC_r(\bG) \ \stackrel{\sim}{\dashrightarrow} \ \bP V_r(\Gr)$
 of Theorem \ref{thm:rho-lambda} for each $r > 0$.  Then the following 
 diagram commutes
  \begin{equation}
 \label{eqn:change-r}
 \begin{xy}*!C\xybox{%
 \xymatrix{\bP\cC_r(\bG)  \ar[r]^-{\rho_r \circ \Lambda_r} \ar[d] & \bP V(\Gr) \ar[d] \ar[r]^-{\Phi \circ \Psi}
 & \Pi(\Gr)  \ar[d]\\
 \bP\cC_{r+1}(\bG)  \ar[r]^-{\rho_{r+1} \circ \Lambda_{r+1}}  & \bP V(\bG_{(r+1)}) \ar[r]^-{\Phi \circ \Psi}
 & \Pi(\bG_{(r+1}) , }
   }\end{xy}
 \end{equation}
 where the  horizontal maps are bicontinuous algebraic maps and the vertical maps are the natural embeddings.
 
 Moreover, for every $\bG$-module $M$ and every $r > 0$, the composition 
 $\Phi \circ \Psi\circ (\rho_r \circ \lambda_r): \bP \cC_r(\bG) \ \to \ \Pi(\bG)$ restricts to a bijection
 \begin{equation}
 \label{eqn:relateM}
 \Phi \circ \Psi\circ (\rho_r \circ \lambda_r): \bP \cC_r(\bG)_M \ \stackrel{\sim}{\to} \  \Pi(\Gr)_{M_{|\Gr}}.
 \end{equation}
 
 Taking colimits with respect to $r$ determines a homeomorphism derived from bicontinuous algebraic map
 \begin{equation}
\label{eqn:cC-iso}
\Phi:  \bP \cC(\bG) \ \stackrel{\sim}{\to} \ \Pi(\bG), \quad \text{restricting to} \quad \bP \cC(\bG)_M 
\ \stackrel{\sim}{\to} \ \Pi(\bG)_M
\end{equation}
for any $\bG$-module $M$.

For each field extension $K/k$, the adjoint action of $\bG(K)$ on the $K$ points of each $\bP\cC_r(\bG)$ 
determines an adjoint action of $\bG(K)$ on the $K$ points $\bP\cC(\bG)$ 
\end{thm}
 
 \begin{proof}
To verify the commutativity of left hand square of (\ref{eqn:change-r}), it suffice to check at geometric points.  This
check is done by comparing $a* \langle \ul B,0 \rangle$ and $\langle a*\ul B, 0 \rangle$
using (\ref{eqn:A1-action}) and (\ref{eqn:A*-action}), where $\langle -,0 \rangle$ sends an $r$-tuple
to the $r+1$-tuple obtained by adding  0 to the right.  The equality (\ref{eqn:epsilons} implies the commutativity
of the right hand square of (\ref{eqn:change-r}).

To verify that $\Phi \circ \Psi\circ (\rho_r \circ \lambda_r): \bP \cC_r(\bG) \ \to \ \Pi(\bG)$ restricts 
as in (\ref{eqn:relateM}), we use the explicit definitions of $\bP \cC(\bG)_M$ in Definition \ref{defn:geom-Cr}
and  of $\Pi(G)_M$ in Recollection \ref{recollect:VrG}.  This restriction is verified by  juxtaposing 
the determination of 
$\Psi \circ (\rho_r \circ \lambda_r): \bP \cC_r(\bG) \to \bP |\Gr|$ in Theorem \ref{thm:rho-lambda}(4), the
definition of $\rho_r$ in   Definition \ref{defn:rho}, and the determination of 
$\Phi \circ \Psi: \bP V_r(\Gr) \to \Pi(\Gr)$ in Theorem \ref{thm:1par,pi}.   Observe that this argument 
remains valid for any $M$ even though we consider $\bP |\Gr|_M$ only for finite dimensional $\Gr$-modules.
 \end{proof}

 To formulate certain functoriality properties (with respect to $\bG$) of $M \mapsto \bP\cC(\bG)_M$, \
 we require the following definition.
 
 \begin{defn}
 \label{defn:functor-exp}
 Let $(\bG, \cE)$ and $(\bG^\prime, \cE^\prime)$ be linear algebraic groups of exponential type.
 A smooth closed embedding $\bG \hookrightarrow \bG^\prime$ is said to be an embedding of 
 exponential type if  $\cE^\prime$ restricted along $f\times id$ to $\cN_p(\bG) \times \bG_a$ 
 equals $\cE: \cN_p(\bG) \times \bG_a \to \bG$.
 \end{defn}

 \begin{prop}
 \label{prop:further-cC} 
 Let $f: (\bG,\cE) \ \to \  (\bG^\prime,\cE^\prime)$ be an embedding of linear algebraic groups of 
 exponential type.
 \begin{enumerate}
 \item
 The embedding $df: \fg \hookrightarrow \fg^\prime$ determines morphisms $\bP \cC_r(\bG) \to \bP \cC_r(\bG^\prime)$
 sending $\ul B$ to \ $df_*(\ul B))$, with
  colimit sending $K$-points of $\bP \cC(\bG)$ to  $K$-points of $\bP \cC(\bG^\prime)$.
    \item
 Composition with $f$ determines $f_*: \Pi(\bG)_{f^*(N)} \ \to \ \Pi(\bG^\prime)_N$ for any $\bG^\prime$-module $N$.
    \item
 Moreover, $f_*(\Pi(\bG)_{f^*(N)})  \ = \ f_*(\Pi(\bG)) \cap \Pi(\bG^\prime)_N$.
    \item
If $J^\prime$ is a mock injective $\bG^\prime$-module, then $f^*(J^\prime)$ is a mock injective $\bG$-module.
\end{enumerate}
 \end{prop}
 
 \begin{proof}
 The first statement is self-evident.  Since  \ $f _*\circ \cE_{\ul B*} \circ \epsilon_r: 
 K\bG_{a(1)} \to K\bG_a \to KG \to KG^\prime$ \ equals
      \ $\cE^\prime_{df_*(\ul B)} \circ \epsilon_r: K\bG_{a(1)} \to K\bG_{a} \ \to \ K\bG^\prime$, 
the action of $(\cE_{\Lambda_r(\ul B})_*(u_{r-1})$ on $f^*N$ equals the
    action on $N$ by the image of $(\cE_{\Lambda_r(\ul B})_*(u_{r-1})$ under $f_*: K\bG \to K\bG^\prime$.
This immediately implies assertions (2) and (3).  Assertion (4) following immediately from assertion (3)
and the detection property of Theorem \ref{thm:properties}(7).
 \end{proof}
 
  \vskip .1in
  
  \begin{defn}
\label{defn:Frob-twist}
We recall that base change along the $p^r$-th power map $k \to k$ associates to a scheme $X$ over
$k$ a map $F^r: X \ \to \ X^{(r)}$.  For a group scheme $G$ and a $G$-module $M$ this leads to 
the definition of the {\it Frobenius twist} $M^{(r)}$ of $M$ given as the restriction along $F^r$ of the
$G^{(r)}$-module $M^{(r)}$. 

If the group scheme $G$ is defined over $\bF_{p^r}$, then we can identify $G^{(r)}$ with $G$ and we 
can associate to a $G$-module $M$ the {\it external Frobenius twist} $M^{[r]}$ defined by 
the restriction along $F^r$ of the $G$-module $M$.  If the action of $G$ on $M$ is defined over $F_{p^r}$,
then $M^{(r)} \ \simeq \ M^{[r]}$ (see \cite[II.3.16]{J}).
\end{defn}

\vskip .1in  
 
 The following proposition  is complementary to Proposition \ref{prop:further-cC}.  One major difference is
 $F^r: \bG \to \bG$  is far from smooth; in fact, its associated differential
 map on Lie algebras is the 0-map.  Another difference is that Proposition \ref{prop:JT-twist} gives an 
 explicit relationship between elements of 
 the support varieties of $M^{[r]}$ and of $M$.

 \begin{prop}
 \label{prop:JT-twist}
 Let $(\bG,\cE)\ \hookrightarrow \ (GL_N,exp)$ be an embedding of exponential type defined 
 over some finite field $\bF_{p^d}$.  Then for any $\bG$-module $M$, 
any $p$-nilpotent element $B$ of $\fg_K$ and any $r \geq d$, the action of $(\cE_{B^{(r)}})_*(u_s)$ on $M_K$
equals the action $(\cE_B)_*(u_{r+s})$ on $M_K^{[r]}$ and the action of $(\cE_B)_*(u_{s})$ on $M_K^{[r]}$
is trivial for $s < r$.

Consequently, for any $r \geq d$, the $K$-points of  $\bP \cC(\bG)_{M^{[r]}}$ are the $K$-points of 
\ $\Proj(\delta_r^{-1}(\cC(\bG)_M))$, where $\delta_r: \cC(\bG) \to \cC(\bG)$ sends a $K$-point 
$(A_0,\ldots,A_{r-1},B_0,\ldots,B_n,\ldots)$ to $(B_0,\ldots,B_n,\ldots)$.
 \end{prop}
 
\begin{proof}
The action $(\cE_B)_*(u_{r+s})$ on $M_K^{[r]}$ is given by the composition
$$M_K \ \to \ M_K \otimes K[\bG] \ \stackrel{1\otimes F^r}{\to} \ M_K \otimes K[\bG] \ \stackrel{(\cE_B)^*}{\to}
\  M_K \otimes K[t] \ \stackrel{1 \otimes u_s}{\to} M_K.$$
Our condition on $(\bG,\cE)$ as the restriction (defined over $\bF_{p^r}$) of the exponential structure 
on $GL_N$ implies that the
composition $F^r \circ \cE_B: \bG_a \to \bG \to \bG$ equals $\cE_{B^{(r)}} \circ F^r: \bG_a \to \bG_a \to \bG$.
(See \cite[Prop1.11]{F2}.)
Thus, the equality of the actions of $(\cE_{B^{(r)}})_*(u_s)$ on $M_K$ and  $(\cE_B)_*(u_{r+s})$ on $M^{[r]}_K$ 
and the triviality of the actions of $(\cE_B)_*(u_{s})$ on $M^{[r]}_K$  for $s < r$ both
follow from the identification of $u_{s+r}$ with $F^r_*(u_s)$.  

This immediately implies the identification of $K$-points of $\cC(\bG)_{M^{[r]}}$ with  $\delta_r^{-1}(\cC(\bG)_M)$
in terms of the $K$-points, and thus the asserted identification of $\bP \cC(\bG)_{M^{[r]}}$.
\end{proof}

The following is an immediate consequence of Proposition \ref{prop:JT-twist} and Theorem \ref{thm:properties}(3),(7).

\begin{cor}
\label{cor:M-twist}
Let $(\bG,\cE)\ \hookrightarrow \ (GL_N,exp)$ be an embedding of exponential type defined 
 over some finite field $\bF_{p^d}$, let
 $r \geq d$, and let $J$ be a mock injective $\bG$-module.  
 \begin{enumerate}
 \item
 $\bP \cC(\bG)_{J^{[r]}} \ \hookrightarrow  \ \bP \cC(\bG)$ equals \ $\bP \cC(\Gr) \ \hookrightarrow  \ \bP \cC(\bG)$.
 \item
 If $M$ is a $\bG$-module, then \ $\bP\cC(\bG)_{M\otimes J^{[r]}}  \ \hookrightarrow  \ \bP\cC(\bG)$
 equals \\
 $(\bP\cC(\bG)_M \cap \bP\cC(\Gr))  \ \hookrightarrow  \ \bP\cC(\bG).$
 \end{enumerate}
\end{cor}

   \vskip .1in

Let $\bG$ be a linear algebraic group of exponential type.  A $\bG$-module $M$ has 
{\it exponential degree $< p^r$} if every  $u_s \in K\bG_a$ with $s \geq r$ acts trivially on 
$\cE_B^*(M)$ for every $B$ a geometric point of  $\cN_p(\fg)$ (see \cite[Defn4.5]{F2}).

\vskip .1in

As we see in Proposition \ref{prop:exp-deg}(1), $\bG$-modules of bounded exponential degree 
have a strong condition on their $\bP \cC_r(\bG)$-support.
Observe that other $\bG$-modules also  satisfy this condition.  For example, if $M$ has 
exponential degree $< p^r$ and if $M \hookrightarrow J$
is an embedding of $M$ in a mock injective $\bG$-module $J$, then $J/M$ also 
has exponential degree $< p^r$ by the ``two out of three" property of Theorem \ref{thm:properties}.4.
   
   \begin{prop}
   \label{prop:exp-deg}
   Let $(\bG, \cE)$ be an algebraic group of exponential type and $M$ be a $\bG$-module.
    \begin{enumerate}
    \item
    If $M$ has exponential degree $< p^r$, then \ $\bP\cC(\bG)_M \ = \ \Proj(\pi_r^{-1}(\pi_r(\cC_r(\bG)_M )))$,
    where $\pi_r: \cC(\bG) \to \cC_r(\bG)$ sends the $K$-point $(B_0,\ldots,B_n,\ldots)$ to $(B_0,\ldots,B_{r-1})$.
   \item
   Any finite dimensional $\bG$-module has bounded exponential degree.
   \item
   The subcategory of $\bG$-modules of exponential degree $< p^r$ is an abelian subcategory of the category
   $Mod(\bG)$ of $\bG$-modules which is closed under extensions and tensor products by
   arbitrary $\bG$-modules.
    \end{enumerate}
    \end{prop}

    \begin{proof}
    Assume that $M$ has exponential degree $< p^r$.  Then 
    assertion (1) follows immediately from the observation for any $t \geq r$ that the action of 
    $\sum_{s=0}^{t-1} \cE_{B_s*}(u_s)$ on $M_K$ equals that of 
    $\sum_{s=0}^{r-1} \cE_{B_s*}(u_s)$ on $M_K$. 
    
    Assertion 2 is proved in \cite[Prop2.6]{F2}.
    
 A straight-forward verification shows that if $M$ has exponential degree $< p^r$ and if 
 $N$ has exponential degree $< p^s$, then any extension of $M$ by $N$ has exponential 
 degree $<  p^{max(r,s)}$ and the tensor product $M\otimes N$ has exponential degreee
 $< p^{r+s}$.
    \end{proof}
    
 Observe that if $\bG$ is as in 
Proposition \ref{prop:JT-twist} and if the $\bG$-module does not have exponential degree $< p^s$,
then $M^{[r]}$ does not have exponential degree $< p^{r+s}$.
       
\vskip .2in


    \section{Examples of $M \ \mapsto \ \bP \cC(\bG)_M$}
    \label{sec:Gexamples}

    In the examples below, we frequently encounter projectivized versions
    of the projection $\pi_r: \cC(\bG) \ \to \ \cC(\Gr)$ sending a $K$-point 
   $(B_0,\ldots,B_n,\ldots)$ to  $(B_0,\ldots,B_{r-1})$.   Starting with
    a subspace $X \subset \bP\cC_r(\bG))$, we abuse notation  by writing  \ 
   $W \ = \ \pi_r^{-1}(X) \ \subset \ \bP\cC(\bG)$  
    determined by the condition that a $K$-point $\ul B$ of $\cC(\bG)$
    represents a $K$-point of $W$ if and only if 
    $(B_0,\ldots,B_{r-1})$ represents a $K$-point of $X$.
    For examples, see Proposotion \ref{prop:exp-deg}(1).

    \begin{ex}
    \label{ex:Ga}
    We consider the special case $\bG = \bG_a$.
    The variety $\cC(\bG_a)$ is identified with $ \bA^\infty$.  The 1-parameter subgroup 
    $\cE_{\ul b}: \bG_{a,K} \to \bG_{a,K}$
    given by the $K$-point $\ul b = (b_0,b_1,\ldots b_n.\ldots )$ of $\bA^\infty$  
    is defined by   $(\cE_{\ul b})^*: K[\bG_a] \  \to \ K[\bG_a]$ sending $T \in K[T] = K[\bG_a]$
    to $b_0T + b_1T^p \ldots b_nT^{p^n} + \ldots $; the $\bA^1$-action is given by sending a 
    $K$-point $(a,\ul b)$ to $(a\cdot b_0,a^p\cdot b_1,\ldots,a^{p^n}\cdot  b_n\ldots)$.
    For any $\bG_{a(r)}$-module $M$, the restriction along the quotient map of algebras 
    (but not a map of Hopf algebras) $k\bG_a \twoheadrightarrow 
    k\bG_{a(r)}$ provides $M$ with the structure of a $\bG_a$-module such that 
     $\cC(\bG_a)_M \ = \ \pi_r^{-1}(\cC_r(\bG_a)_M) \ \subset \ \cC(\bG_a)$
     (where $\pi_r: \cC(\bG_a) \to \cC_r(\bG_a)$ is the projection onto the first $r$ factors).
    If $ W \ \subset \ \bP^{r-1} = \bP \cC_r(\bG_a)$ is any closed subset, some tensor product of 
    Carlson $L_\zeta$-modules for $\bG_{a(r)}$, $M =  \otimes_i L_{\zeta_i}$,
    has the property that $\Pi(\Gr)_M = W$ and  thus\ $\bP \cC(\bG_a)_M \ = \pi_r^{-1}(W)$.  
    Consequently, a closed subset $W\ \subset \ \bP\cC(\bG_a)$
    is of the form $\bP\cC(\bG_a)_M$ for some finite dimensional $\bG$-module $M$
     if and only if $W \ = \ \Proj \pi_r^{-1}(\pi_r(W))$  for some $r>0$.
        
    More generally, using Rickard idempotent modules, we may realize any subset $X \ \subset \ \bP \cC_r(\bG_a)$
    as $\Pi(\bG_{a(r)})_M$ for some (not necessarily finite dimensional) $\bG_{a(r)}$-module $M$.  Thus, given 
    any subset $X \ \subset \ \bP \cC_r(\bG_a)$, we may find some $\bG_a$-module $M$ with the property
    that $\bP \cC(\bG_a)_M \ = \ \pi_r^{-1}(X)$.
    
    We explicitly construct infinite dimensional $\bG_a$-modules $M$  with the property that $X = \bP \cC(\bG_a)_M$
    does not satisfy \ $X = \pi_r^{-1}(W)$ for any $W \subset \bP\cC_r(\bG_a)$.
    Let $L$ denote the regular representation of $\bG_a$
    on $k[\bG_a] \simeq k[T]$.  Thus, $k\bG_a = k[u_0,\ldots,u_n,\ldots]/(u_i^p)$ acts on $k[T]$ with each $u_i$
    a derivation and $u_i(T^{p^j}) = \delta_{i,j}$.  
    For any subset $S \subset \bN$, we define the $\bG_a$-module $L_S$
    by identifying $L_S$ with $k[T]$ as a $k$-vector space and defining  $u_i$ on $L_S$ to be 0 if $i \notin S$
    and $u_i$ as the derivation $u_i(T^{p^j}) = \delta_{i,j}$ if $i \in S$.  For any $\ul B \not= 0$, $\cE_{\ul b}$ 
    represents a  
    $K$-point of $\bP \cC(\bG_a)_{L_S}$ if and only if $b_i = 0, \ \forall i \in S$.  This is verified by observing
    that if  $b_i = 0$ for all $ i \in S$, then the Jordan type of $L_S$ at $\cE_{\ul b}$ consists only of blocks of size 1, 
    whereas if  $b_i \not= 0 $ for some $i \in S$, then $\E_{\ul b}$ determines a 
    non-zero $\pi$-point of $\bP\cC(\bG_{a(i+1})_{L_S}$.
    
    The preceding analysis applies with little change if one replaces $\bG = \bG_a$ by $\bG = \bG_a^{\times s}$
    for any positive integer $s$.
 \end{ex}
    
    \vskip .1in
    
     \begin{ex}
    \label{ex:St-s}
  Let $\bG$ be a split, reductive algebraic group which admits an embedding 
  $(\bG,\cE) \ \hookrightarrow \ (GL_N,exp)$  
  of exponential type.  For any $s > 0$, let $St_s$ be the $\bG$-module given by 
  $$St_s \quad \equiv \quad L((p^s-1)\rho) \quad = \quad H^0(\bG/\bB,(p^s-1)\rho)).$$
  The restriction of $St_s$ to $\bG_{(s)}$ is both irreducible and injective.
  Consequently, the only $K$-point of $\cC_s(\bG)_{St_s}$ is 0.  
  Since $St_s$ is a block of the $\bG$-module $k\bG_{(s)}$, $u_r$ acts trivially on
  $St_s$ for $r \geq  s$ .
  Thus, \ $\bP \cC(\bG)_{St_s} \ = \ \pi_s^{-1}(\{ \emptyset \})$; in other words, $\bP \cC(\bG)_{St_s}$
  is the center of the projection $\pi_s: \bP\cC(\bG) \dashrightarrow \bP \cC_s(\bG)$.
  
  Following \cite[Prop7.8]{SFB2}, we can use this computation to determine
  \ $\bP\cC(\bG)_{H^0(\lambda)}$ \ for any dominant weight $\lambda$ of the form $n\rho$.
  \end{ex}
    
    \vskip .1in

    \begin{ex}
    \label{ex:U3}
    Let $\bU_3$ denote the Heisenberg group, the unipotent radical of a split Borel 
    subgroup of $GL_3$.   The computations in \cite{F5} give the 
    identification  \ $ \cC_r(\bU_3) \ = \ Y_r \times \bA^r \ \subset \ \bA^{2r} \times \bA^1, $ \
    where $Y_r$ is the closed subvariety of $\bA^{2r}$ with generators 
    $x_{1,2}(\ell),\ x_{2,3}(\ell^\prime), 0 \leq \ell, \ell^\prime < r$ subject to
    the relations 
    \ $x_{1,2}(\ell)\cdot x_{2,3}(\ell^\prime) \ = \ x_{2,3}(\ell)\cdot x_{1,2}(\ell^\prime), \quad 0 \leq \ell < \ell^\prime.$ \
    In fact, $Y_r$ is shown to be a rational variety of dimension $r+1$, smooth outside of the origin \cite[Prop5.2]{F5}.
    
    Thus, \ $\bP \cC_r(\bU_3) \ = \ (\bP Y_r) \# \bP^{r-1} \ \subset \ \bP^{3r-1}$.  Furthermore, the colimit 
    $\bP \cC(\bU_3)$ is identified with \ $(\bP Y_\infty)  \# \bP^\infty$.    
    
    Observe that the restriction to $\cC_r(\bU_3) $ of the adjoint 
    action of $\bU_3$ to $\cN_p(\fu_3)^{\times r}$ is non-trivial.
    Since the adjoint action must stabilize subsets of the form $\bP\cC(\bU_3)_M$ for a $\bU_3$-module $M$
    by Theorem \ref{thm:rho-lambda}(3), this constrains which subspaces are of the form $\bP\cC(\bU_3)_M$.

    One class of examples of $\bU_3$-modules is given by inflation $pr^*(N)$ of 
    $\bG_a^{\times 2}$-modules 
    $N$ along the projection homomorphism $pr: \bU_3 \to \bG_a^{\times 2}$.  This projection 
    induces the projection map 
  $$p_r: \bP\cC_r(\bU_3) \ \simeq \ (\bP Y_r) \# \bP^{r-1} \ \dashrightarrow \ \bP Y_r   \ \hookrightarrow \bP^{2r-1}.$$
    Consequently, for any (Zariski) closed subspace $W \subset \bP Y_r$, we may realize
    $p_r^{-1}(W) \subset \bP\cC_r(\bU_3)$ as $\bP\cC(\bU_3)_{p_r^*(N)}$ for some finite dimensional
    $\bU_3/Z$-module $N$ by Example \ref{ex:Ga}.
    
    Interesting infinite dimensional examples are given by induction from the 
    central subgroup $\bG_a \ \simeq \ Z \ \subset \ \bU_3$.
    For example, let $M = \ind_Z^{\bU_3} k \simeq k[\bU_3/Z]$.  The action of 
    $\bU_3$ on $M$ is given by the projection to $\bU_3/Z$ followed by the 
    left regular representation of $\bU_2/Z$.
    Thus, if $\cE_{\ul B} \in \cC_r(\bU_3)$ factors through $Z$, then the action of 
    $(\cE_{\ul B})_*(u_{r-1})$ on $M$ is trivial and 
    therefore $\ul B \in \cC_r(\bU_3)$ lies in $\cC_r(\bU_3)_M$.
    Otherwise, the composition $pr \circ \cE_{\ul B}: \bG_a \ \to \  \bU_3 \ \to \ \bU_3/Z$ 
    represents a non-zero point of $\cC_r(\bU_3/Z)$; since $\cC(\bU_3/Z)_{k[\bU_3/Z]} = 0$,
    the action at the 1-parameter subgroup $\cE_{\ul B} \to \bU_3$ on $M$ can be 
    identifiied with the action at $p_r \circ \cE_{\ul B} \to \bU_3/Z$ for the regular 
    representation of $\bU_3/Z$ which is free.
    We conclude that \ $Z \hookrightarrow \bU_3$ \ induces a homemorphism
    \ $\bP\cC(Z) \ \stackrel{\sim}{\to} \ \bP\cC(\bU_3)_M$.
        \end{ex}
    
    \vskip .1in
    
    \begin{ex}
    \label{ex:reductive}
    Let $\bG$ be a simple algebraic group over an algebraically closed field $k$ 
    with Coexter number $h< p$.  Assuming the  Lusztig
    character formula is valid for all restricted dominant weights $\lambda$, C. Drupieski, D. Nakano, and B. Parshall 
    showed in \cite[Thm3.3]{DNP} for any restricted dominant weight $\lambda$ that
    $$\cC_1(\bG)_{L(\lambda)} \ = \ \bG \cdot \fu_{J(\lambda)}$$
    where $\fu_{J(\lambda)}$ is the Lie algebra of the unipotent radical of a parabolic subgroup $\bP_{J(\lambda)}$
    explicitly determined by $\lambda$ and $L(\lambda)$ is the irreducible $\bG$-module of 
    high weight $\lambda$. 
    
    For $\bG$ a classical simple algebraic group of rank $\ell$ satisfying $p > \frac{\ell^2}{2} + 1$,
 P. Sobaje in \cite{Sobj} explicitly computed $\cC_r(\bG)_{L(\lambda)}$
 for any dominant weight $\lambda = \sum_{i \geq 0}p^i\lambda_i$ (with each $\lambda_i$ 
 a restricted dominant weight).  The form of 
 Sobaje's determination is 
    $$\bP\cC_r(\bG)_{L_\lambda} \ = \ \{ \ul B \in \bP\cC_r(\bG): B_i \in \bG\cdot \fu_{J(\lambda_i)} \}.$$
    
   Assume further that $\lambda = \sum p^i\lambda_i$ satisfies the condition \ 
    $2\sum_{j=1}^\ell \langle \lambda_i,\omega_j^\vee \rangle < p(p-1)$ \ for each $i$.
Then, as argued in \cite[Prop5.1]{F2}, we can identify the $K$-points of $\cC(\bG)_{L_\lambda}$
as those $K$-points $\cE_{\ul B}$ of $\cC(\bG)$ such that $B_i$ is a $K$-point  of 
$\bG\cdot \fu_{J(\lambda_i)}.$
 
 Similar results are obtained by these authors with the irreducible module $L(\lambda)$
 replaced by $H^0(\lambda) \ = \ ind_{\bB}^{\bG} (k_\lambda)$, where $\bB \ \subset \bG$ is 
 a Borel subgroup of $\bG$.   
    \end{ex}

\vskip .2in


    \section{Stable module categories}
\label{sec:stab-cat}

    For a finite group scheme $G$, a natural domain for support theory is the stable module
    category $StMod(G)$, with triangulated category structure introduced by Happel \cite{Ha}.  
    This triangulated structure utilizes the fact that the group algebra $kG$ is self-injective, implying that projective
    objects are the same as injective objects in the abelian category $Mod(G)$.  
Our primary objective
    in this section is to introduce in Definition \ref{defn:stmod} the triangulated categories
     $StMod(\bG)$ and $stmod(\bG)$ associated to the triangulated structures of the homotopy 
     categories of (cochain) complexes $\cK^b(Mod(\bG))$ and $\cK^b(mod(\bG))$.

  \vskip .1in
  
  We begin by recalling various constructions associated to the abelian category $Mod(G)$
 of $G$-modules and the full abelian subcategory $mod(G)$
 of finite dimensional $G$-modules for an affine group scheme $G$ over $k$.
 We denote by $CH(Mod(G))$ the abelian category of (cochain) complexes of $G$-modules
and by  $\cK(Mod(G))$ the associated homotopy category with respect to cochain homotopy equivalence.
The additive category $\cK(Mod(G))$ has the natural structure of a triangulated category with 
basic exact triangles $C^\bu \to D^\bu \to E^\bu \to C^\bu[1]$ associated to short 
exact sequences $0 \to C^\bu \to D^\bu \to E^\bu \to 0$ in $CH(Mod(G))$.  We are primarily
interested in the analogous structures $CH^b(Mod(G)), \ \cK^b(Mod(G))$ whose objects are bounded 
complexes of $\bG$-modules and their full subcategories $CH^b(mod(G)), \ \cK^b(mod(G))$
of bounded complexes of finite dimensional $G$-modules,

We say that a bounded complex $C^\bu$ of $G$-modules in $CH^b(Mod(G))$
has length  $\leq d$ if there exists some integer $m$  such that $C^i = 0$ for $i < m$ or $i > m+d$.
 In particular, complexes $M[n]$ which have the $G$-module $M$ in some degree $n$ and 
 0 in all other degrees have length 0.
 
We recall that the stable module category $StMod(G)$ of a finite group scheme $G$ is the
additive category whose objects are $G$-modules and whose abelian group of maps 
$Hom_{StMod(G)}(M,N)$ is the quotient of $Hom_{Mod(G)}(M,N)$ by the subgroup of maps 
$M \to N$ which factor through an injective $G$-module.  One equips $StMod(G)$ with
a triangulated structure whose +1 shift is given by sending $M$ to the
cokernel $\Omega^{-1}(M)$ of a minimal embedding $M \hookrightarrow I_M$ with $I_M$
an injective $G$-module, and whose -1 shift is given by sending $M$ to the kernel $\Omega^1(M)$
of a minimal surjection $P_M \twoheadrightarrow M$ with $P_M$ a projective $G$-module.  
We define $stmod(G) \subset StMod(G)$
to be the tensor triangulated subcategory whose objects are finite dimensional $G$-modules.

If $\cC \ \hookrightarrow \ \cD$ is an embedding of triangulated categories (which implies by
definition that the image of $\cC$ is a full subcategory of $\cD$), then the Verdier
quotient $\cD \ \to \  \cD/\cC$ is the universal triangulated map from $\cD$ to
a triangulated category $\cE$ with kernel containing $\cC$.  If $\cC$ is a thick subcategory
(i.e., containing any summand in $\cD$ of an object of $\cC$), then the kernel of 
$\cD \to \cD/\cC$ equals $\cC$.  The Verdier quotient is constructed as 
the category whose objects are the same as objects of $\cD$ and whose maps $X \to Y$ 
in $\cD/\cC$ are
equivalence classes of ``roofs": triples $X \leftarrow E \to Y$ in $\cD$ such that the cone of $E \to X$
is an object of $\cC$.   We refer the reader
to \cite{Nee} for a detailed discussion of Verdier quotients of triangulated categories.

	Granted the construction of the Verdier quotient, we readily observe that  $\cD \to \cD/\cC$ 
is the localization $\cD \to \cD[mor_{\cC}^{-1}]$,  where $mor_\cC$ is the class of maps $X \to Y$ 
in $\cD$ whose cone lies in $\cC$.  Namely, $mor_\cC$ is a ``multiplicative system"
as considered for example in \cite[10.3.4]{Wei}.  Since $\cD \to \cD/\cC$ sends each $X \to Y$
 in $mor_\cC$ to an invertible map,  $\cD \to \cD[mor_{\cC}^{-1}]$ is well defined and factors
$\cD \ \to \ \cD/\cC$; on the other hand, granted the existence of $\cD[mor_{\cC}^{-1}]$,
the fact that $\cC$ lies in the kernel of $\cD \to \cD[mor_{\cC}^{-1}]$ provides the inverse functor
$\cD/\cC \to \cD[mor_{\cC}^{-1}]$ factoring $\cD \to \cD[mor_{\cC}^{-1}]$.

	We remind the reader that the bounded derived category of an abelian category $\cA$, $D^b(\cA)$,
is the Verdier quotient of the homotopy category $\cK^b(\cA)$ of bounded chain complexes of $\cA$ 
by the thick triangulated subcategory $Quasi^b(\cA)$ of bounded complexes quasi-isomorphic to 0
(i.e., which are exact).   If $\cA$ has enough injective objects (for example, if $\cA$ is the category
of $G$-modules for an affine group scheme), then $D^b(\cA)$ can be defined more simply by first
explicitly defining the derived category $D^+(\cA)$ of complexes bounded below using injective 
resolutions and then restricting to objects which are bounded complexes.  
If  $f: I^\bu \to J^\bu$ is a quasi-isomorphism (i.e., induces an isomorphism in
cohomology) of bounded complexes of injective $G$-modules, then $f$ 
is a chain homotopy equivalence.   Consequently, the homotopy category $\cK^\bu(Inj)$ of bounded complexes for
the full additive subcategory $Inj(Mod(G)) \subset Mod(G)$ embeds as a full subcategory of $D(Mod(G))$.

    The following well-known theorem (see, for example, \cite[Thm2.1]{R1}) provides motivation
for our constructions.  We recall that a perfect complex of $G$-modules for a finite group
scheme $G$ is a bounded complex of finite dimensional $G$-modules quasi-isomorphic to a 
bounded complex of finite dimensional injective $G$-modules.   Thus, the subcategory 
$Perf(G) \ \subset \ D^b(mod(G))$ of perfect complexes is the essential image in $D^b(mod(G))$
of the homotopy category of bounded complexes of finite dimensional injective $G$-modules.
$\cK^b(inj(G)) \ \subset \ D^b(mod(G))$.

    \begin{thm} 
    \label{thm:rickard}
If $G$ is any finite group scheme, then the natural map 
    $mod(G) \ \to \ \cK^b(mod(G))$ induces an equivalence of tensor triangulated categories
    $$\Psi: stmod(G) \quad \stackrel{\sim}{\to} \quad D^b(mod(G))/Perf(G).$$
  \end{thm}

This theorem is somewhat surprising, even in the case of $\bG_{a(1)}$ (i.e., $k[t]/t^p$-modules).
Consider the complex $C^\bu$ with $C^0 = k[t]/t^p$, $C^1 = k[t]/t^p$ and 
$d: C^0 \to C^1$ given by multiplication by $t$;  thus $C^\bu \not= 0$ in $D^b(mod(G))$
but is 0 in $stmod(G)$.   In the proof of  Theorem \ref{thm:St-rickard}, the construction of 
 $\tilde \Phi: D^b(mod(G)) \to stmod(G)$ applied to the complex $C^\bu$
sends $C^\bu$ to the cone of $\tilde \Phi(C^0) = 0 \in stmod(G)$ mapping to  
$\tilde \Phi(C^1) = 0  \in stmod(G)$ which is 0.
\vskip .1in

We shall utilize the truncation functors
$$\tau_{\leq n}, \tau_{\geq m}: CH(Mod(G)) \ \to CH(Mod(G)), \quad 
C^\bu \mapsto \tau_{\leq n}(C^\bu) \hookrightarrow C^\bu, \quad C^\bu \twoheadrightarrow \tau_{\geq m}(C^\bu).$$
Here, $\tau_{\leq n}(C^\bu)$ is the subcomplex of $C^\bu$ such that $\tau_{\leq n}(C^\bu)^i = 0$ 
for $i > n$, $\tau_{\leq n}(C^\bu)^n  = ker\{d^n\}$, and $\tau_{\leq n}(C^\bu)^i  = C^i$ for $i < n$;
and $\tau_{\geq m}(C^\bu)$ is the quotient complex of $C^\bu$ such that $\tau_{\geq m}(C^\bu)^n = 0$ 
for $n < m$, $\tau_{\geq m}(C^\bu)^m = C^m/im\{ d^{m-1}\}$, and $\tau_{\geq m}(C^\bu)^n = C^n$ for $n > m$.
So defined, if $C^\bu \to D^\bu$ is a quasi-isomorphism, then $\tau_{\leq n}(C^\bu) 
\to \tau_{\leq n}(D^\bu)$ and $\tau_{\geq m}(C^\bu) \to \tau_{\geq m}(D^\bu)$ are also quasi-isomorphisms.

In the following theorem, we consider the thick tensor triangulated subcategory
$\cI nj^b(Mod(G)) \ \subset \ D^b(Mod(G))$ of complexes quasi-isomorphic to bounded complexes
of injective $G$-modules.   This is the ``essential image" of $\cK^b(Inj(G)) \ \subset D^b(Mod(G))$.

\begin{thm}
\label{thm:St-rickard} (see \cite[Thm2.1]{R1})
For any finite group scheme $G$, the natural map 
$Mod(G) \ \to  \ \cK^b(Mod(G))$ (sending the $G$-module $M$ to the 
chain complex $M[0]$ concentrated in degree 0) induces an equivalence of 
tensor triangulated categories
$$\Psi: StMod(G) \quad \stackrel{\sim}{\to}   \quad D^b(Mod(G))/\cI nj^b(Mod(G)).$$

Furthermore, we construct $\tilde \Phi: D^b(Mod(G)) \ \to \ StMod(G)$ inducing the inverse
\begin{equation}
\label{eqn:PHI}
\Phi = \Psi^{-1}:  D^b(Mod(G))/\cI nj^b(Mod(G)) \quad \stackrel{\sim}{\to} \quad StMod(G)
\end{equation}
sending $M[n]$ to $\Omega^{-n}(M)$.
\end{thm}

\begin{proof}
Necessarily, $\Psi(\Omega^{n}(M)) \ = \ M[-n]  \in D^b(Mod(G))/Inj^b(Mod(G))$.
As explained in Rickard's proof of \cite[Thm2.1]{R1}, an exact triangle  \ 
$X \to Y \to Z \to X[1]$ \ in $StMod(G)$ arises from a pushout  diagram in 
$Mod(\bG_{(r)})$ from the short exact sequence $0 \to X \to I \to X[1] \to 0$ 
(with $I$ injective) along $X \to Y$ to the short exact sequence
$0\to Y \to Z \to X[1] \to 0$.  Consequently, the result of applying $\Psi$ to a such 
an exact triangle in $StMod(G)$,  \ $\Psi(X) \to \Psi(Y) \to \Psi(Z) \to \Psi(X)[1]$,
is isomorphic in $D^b(Mod(G))/Inj^b(Mod(G))$ to the image of the exact 
triangle in $\cK^b(Mod(G))$ arising from the short exact sequence
$0 \to X \to Y\oplus I \to( Y\oplus I)/X \to 0$ in $Mod(G)$.    Thus, $\Psi$ is exact; in 
other words, it sends exact triangles to exact triangles.  

To show that $\Psi$ is a full functor, we first show that $Mod(G) \to D^b(Mod(G))$ is
full, then observe that $D^b(Mod(G)) \to D^b(Mod(G))/\cI nj^b(Mod(G))$ is full, so
that $\Psi$ must also be full.  Let $M, \ N$ be $G$-modules and consider a map
from $M[0]$ to $N[0]$ in $D^b(Mod(G))$ represented by a ``roof" $M[0] \ \leftarrow 
C^\bu \ \to \ N[0]$ in $\cK^b(Mod(G))$.    This roof factors through the quotient map
$C^\bu \twoheadrightarrow \tau_{\geq 0}(C^\bu)$ (which is  a quasi-isomorphism). 
On the other hand, $M \ \simeq \ ker\{ d^0\}/im\{d^{-1}\} \hookrightarrow \tau_{\geq 0}(C^\bu)^0$.
Thus,  there is a natural map  of roofs from 
$(M[0] \stackrel{=}{\leftarrow} M[0] \to N[0)]$ to $(M[0] \ \leftarrow \tau_{\geq 0}(C^\bu) \ \to \ N[0])$, 
thereby establishing the fullness of $Mod(G) \to D^b(Mod(G))$.

We next show that $Hom_{D^b(Mod(G))}(M,N)  \ \to \ Hom_{D^b(Mod(G))/\cI nj^b(Mod(G))}(M,N)$ is 
surjective.   Consider a roof of maps $M[0] \ \leftarrow C^\bu \ \to \ N[0]$ in $D^b(Mod(G))$, with
associated exact triangle $J^\bu \to C^\bu \to M[0]$ in $D^b(Mod(\bG))$ for some bounded
complex $J^\bu$ of injective $G$-modules.   Replacing $C^\bu$ by the mapping 
cylinder of $J^\bu \to C^\bu$ if necessary, we may assume that $J^\bu \hookrightarrow C^\bu$
so that there is an isomorphism $C^\bu \to ((C^\bu/J^\bu) \oplus J^\bu)$ and $C^\bu \to M[0]$ 
factors as $C^\bu \to ((C^\bu/J^\bu) \oplus J^\bu) \to M[0]$ with $C^\bu/J^\bu \to M[0]$ an
isomorphism in $D^b(Mod(G))$.  Thus, the given roof is equal as a map in $D^b(Mod(G))/\cI nj^b(Mod(G))$
to a roof of the type considered in the previous paragraph and thus is the image under 
$\Psi$ of a map $M \to N$.

As shown in the proof of \cite[Thm2.1]{R1}, exactness and fullness of $\Psi$ imply that $\Psi$
is also faithful.  Since any summand of an injective $\bG$-module is itself injective,
we conclude that $\cI nj^b(Mod(G))$ is a thick subcategory of $D^b(Mod(G))$.
To check that every object $C^\bu \in 
\cK^b(Mod(G))/\cI nj^b(Mod(G))$ lies in the image of $\Psi$ (up to isomorphism), 
we argue by induction on the length $d$ of $C^\bu$, recognizing that any 
complex of length $0$ is of the form $M[n]$ and thus is in the image of $\Psi$.
We construct
$\tilde \Phi: D^b(Mod(G)) \ \to \ StMod(G)$ (which induces $\Phi$ of (\ref{eqn:PHI}))
with the property that $\Psi \circ \tilde \Phi$ sends $X^\bu \in D^b(Mod(G))$
to its image in $D^b(Mod(G))/\cI nj^b(Mod(G))$ under the quotient map.  We proceed by induction, assuming that 
$\tilde \Phi$ has been defined on the full subcategory of $D^b(Mod(G))$ whose objects are 
quasi-isomorphic to complexes of length $\leq d$.

Let $X^\bu \in D^b(Mod(G))$ be a complex of length $d+1$ with $X^{n+1} \not= 0$ and 
$X^i = 0$ for $i > n+1$.  We consider the exact triangle in $\cK^b(Mod(G))$
\begin{equation}
\label{eqn:tri1}
H^{n+1}(X^\bu)[n] \to \tau_{\leq n}(X^\bu) \to X^\bu \to H^{n+1}(X^\bu)[n+1],
\end{equation}
and recognize that $\tau_{\leq n}(X^\bu)$ is a complex of
length $\leq d$.   We extend $\tilde \Phi(H^{n+1}(X^\bu)[n]) \to \tilde \Phi( \tau_{\leq n}(X^\bu))$
to an exact triangle
in $StMod(G)$
\begin{equation}
\label{eqn:tri2}
\tilde \Phi(H^{n+1}(X^\bu)[n])  \ \to \ \tilde\Phi(\tau_{\leq n}(X^\bu)) \ \to M \ \to \ \tilde\Phi(H^{n+1}(X^\bu)[n+1]).
\end{equation}
Then $\Psi$ defines an isomorphism of exact
triangles in $D^b(Mod(G))$ from $\Psi$(\ref{eqn:tri2}) to (\ref{eqn:tri1}).
We define $\tilde \Phi(X^\bu) \equiv M$.  In particular, $X^\bu \ = \ \Psi(M)$ is in the 
image of $\Psi$.

Assume we have defined $\tilde \Phi(X^\bu)$ for all complexes of length $d+1$ (with
arbitrary $n$ as top non-vanishing degree).  Let $Y^\bu$ be another bounded 
complex of length $d+1$.
For any map $f: X^\bu \to Y^\bu \in D^b(Mod(G))$, we define 
$\tilde\Phi(f): \tilde\Phi(X^\bu) \to \tilde\Phi(Y^\bu)$ to be the unique map 
satisfying the condition that $\Psi (\tilde \Phi(f))$
equals the image of $f$ in $D^b(Mod(G))/Inj^b(Mod(G))$. 
One readily verifies that $\tilde \Phi$ so defined is a  functor which preserves
exact triangles and induces $\Phi$ inverse to $\Psi$.
\end{proof}

\vskip .1in

\begin{remark}
\label{rem:warning}
Let $G$ be a finite group scheme.  Although the category $StMod(G)$ has arbitrary direct sums, 
the direct sum $\bigoplus_{j\in J} C^\bu_j \ \in \ \cK(Mod(G))$ of an infinite family of bounded 
complexes of $G$-modules does {\it not} represent the 
direct sum in $StMod(G)$.

Arbitrary direct sums in $D^b(G)/\cI nj(Mod(G))$ are determined by using the triangulated
structure to represent each bounded complex $C^\bu_j$  as $M_j[0]$ for some
$\Gr$-module $M_j$ and then taking the arbitrary
sum $\bigoplus_j M_j$ of $G$-modules.
\end{remark}

\vskip .1in
 
 \begin{remark}
 \label{rem:alternate-stable}
 The tensor triangulated category $D^b(Mod(G))/\cI nj^b(Mod(G))$ of Theorem \ref{thm:St-rickard}
 has the following two alternate formulations.  (Analogous formulations are valid for 
 $D^b(mod(G))/Perf(G)$.)
 
 The first is as the Verdier quotient $\cK^b(Mod(G))/\cI nj^b(Mod(G))$, provided that 
 $\cI nj^b(Mod(G))$ is viewed as a thick tensor subcategory of $\cK^b(Mod(G))$ rather
 than as a subcategory of $D^b(Mod(G))$.  This is essentially immediate from the 
 universal property satisfied by a Verdier quotient.
 
 The second is the somewhat more ``elementary" formulation as the category of left (or right) fractions 
 of $D^b(Mod(G))$ for the multiplicative system $\cS \ \subset Mor(D^b(Mod(G)))$ consisting of 
 morphisms of $f: C^\bu \to D^\bu$ whose cone $c(f)$ is quasi-isomorphic to a bounded complex of 
 injective $G$-modules.  (See \cite[Chap10]{Wei}.)
 \end{remark}

We obtain the following useful corollary of Theorem \ref{thm:St-rickard}.
 
 \begin{cor}
 \label{cor:detect-iso}
 Let $G$ be a finite group scheme.  
 \begin{enumerate}
 \item
A complex $C^\bu \in \cK^b(Mod(G))$
is quasi-isomorphic to a bounded complex of injective $G$-modules if and only
if \ $\tilde\Phi(C^\bu) \ = \  0$ in $St(Mod(G))$.  
\item
For any map $f:C^\bu \to D^\bu$ in $CH^b(Mod(G))$, the map $\Phi(f)$ 
 is an isomorphism in $StMod(G)$ if and only if the cone of $f$ in $\cK^b(Mod(\bG))$
 is an object of $Inj^b(Mod(G))$.
 \item
 If $C^\bu \to D^\bu \to E^\bu \to C^\bu[1]$ is an exact triangle in $D^b(Mod(G))$ and 
 if $C^\bu$ is quasi-isomorphic to a bounded complex of injective $G$-modules,
 then $D^\bu$ is quasi-isomorphic to a bounded complex of injective $G$-modules \ $\iff$ \ $E^\bu$ is
 quasi-isomorphic to a bounded complex of injective $G$-modules.
 \end{enumerate}
 \end{cor}
 
 \begin{proof}
 Assertion (1) is immediate from Theorem \ref{thm:St-rickard} and the fact that the kernel of
 the Verdier quotient $\cD \to \cD/\cC$ is $\cC$ whenever the triangulated subcategory
 $\cC \subset \cD$ is a thick subcategory.  Assertion (2) follows immediately from 
 Theorem \ref{thm:St-rickard} and Remark \ref{rem:alternate-stable}.  
 
 Since $\tilde \Phi$ is exact, we conclude that 
 $\tilde \Phi(C^\bu) \to \tilde \Phi(D^\bu) \to \tilde \Phi(E^\bu) \to \tilde \Phi(C^\bu)[1]$ is 
 an exact triangle in $St(Mod(G))$.  Thus, if $\tilde \Phi(C^\bu) = 0$, then 
 $\tilde \Phi(D^\bu) = 0$ $\iff$ $\tilde \Phi(E^\bu) = 0$.  This, together with (1) implies
 assertion (3).
 \end{proof}

In our approach to stable module categories for $\bG$-modules,  mock injective
$\bG$-modules as in Definition \ref{defn:mock} play the role that injective $\Gr$-modules play in defining $StMod(\Gr)$.
   
  \begin{defn}
  \label{defn:cMock}
   Let $\bG$ be a connected linear algebraic group.   We define $\cM ock^b(\bG)$
  to be the full subcategory of $D^b(Mod(\bG))$ consisting
 of those bounded complexes $C^\bu$ of $\bG$-modules each of whose
 restrictions $C^\bu_{| \Gr}$ is an object of $\cI nj^b(\Gr)$.
 
We define $mock^b(\bG)$
  to be the full subcategory of $D^b(mod(\bG))$ consisting
 of those of bounded complexes $C^\bu$ of finite dimensional $\bG$-modules each of whose
 restrictions $C^\bu_{| \Gr}$ is an object of $\cI nj^b(\Gr)$. 
 \end{defn}
 
 \vskip .1in
 
 We easily verify that  $\cM ock^b(\bG)\  \hookrightarrow \ D^b(\bG)$ is a thick tensor ideal by 
 observing that $\cM ock^b(\Gr) \hookrightarrow D^b(\Gr)$ is a thick tensor ideal
 for every $r > 0$.

\begin{defn}
\label{defn:stmod}
 Let $\bG$ be a connected linear algebraic group.   
 We define $StMod(\bG)$ to be the tensor triangulated category defined 
 as the Verdier quotient
 $$StMod(\bG) \quad \equiv \quad D^b(Mod(\bG))/\cM ock^b(\bG)$$
 and denote by \ $q_{Mock}: D^b(Mod(\bG)) \to StMod(\bG)$ the
 quotient functor.
 
 We define $stmod(\bG)$ to be the tensor triangulated category defined as
the Verdier quotient
 $$stmod(\bG) \quad \equiv \quad D^b(mod(\bG))/mock^b(\bG)$$
 and denote by \ $q_{mock}: D^b(mod(\bG)) \to stmod(\bG)$ the
 quotient functor.
 \end{defn}
\vskip .1in

\begin{remark}
\label{rem:inf-inj}
Let $\bG$ be a connected linear algebraic group with the property that 
the minimal dimension of a non-trivial injective $\Gr$-module becomes arbitrarily large 
as $r$ increases.  Then every mock injective $\bG$-module is infinite dimensional, so that
$stmod(\bG) \ = \ D^b(mod(\bG))$.
This is the case for each of the examples in Section \ref{sec:Gexamples} assuming that
$\bG$ does not have a factor which is a torus.
\end{remark}

We shall frequently abuse notation in what follows by denoting by $M$ both a $\bG$-module 
$M$ and its restriction $M_{|\Gr}$ to a Frobenius kernel $\Gr \hookrightarrow \bG$.

\begin{prop}
\label{prop:detect-faithful}
Let $\bG$ be an irreducible linear algebraic group and $ C^\bu, \ D^\bu$ be two 
 bounded complexes of finite dimensional $\bG$-modules.  Then for $r$ sufficiently 
large, the restriction map  is an isomorphism:
\begin{equation}
\label{eqn:bij-st}
Hom_{stmod(\bG)}(C^\bu,D^\bu) \ \to \ Hom_{stmod(\Gr)}(C^\bu_{|\Gr},D^\bu_{|\Gr}).
\end{equation}
\end{prop}

\begin{proof}
For any finite dimensional $\bG$-modules $M, \ N$, the restriction map \\
\ $Hom_{Mod(\bG)}(M,N) \ \to \ Hom_{mod(\Gr)}(M_{|\Gr},N_{|\Gr})$ is an isomorphism for $r > > 0$
by \cite[I.9.8(6)]{J}.  This readily implies that 
\begin{equation}
\label{eqn:Hom-bij-cx}
Hom_{\cK^b(mod(\bG))}(C^\bu,D^\bu) \ \stackrel{\sim}{\to}  \ Hom_{\cK^b(mod(\Gr))}(C^\bu_{|\Gr},D^\bu_{|\Gr})
\quad r >> 0
\end{equation} 
for any pair $C^\bu, \ D^\bu$ of bounded complexes of finite dimensional $\bG$-modules.  
The isomorphism (\ref{eqn:Hom-bij-cx})
clearly preserves quasi-isomorphisms, thus easily implies that the restriction map \ 
\begin{equation}
\label{eqn:Hom-Db-cx}
Hom_{D^b(mod(\bG))}(C^\bu,D^\bu) \ \to \ Hom_{D^b(mod(\Gr))}(C^\bu_{|\Gr},D^\bu_{|\Gr})
\end{equation}
 is also an isomorphism.

Two maps $f,g \in Hom_{D^b(mod(\bG))}(C^\bu,D^\bu)$ are equal in $Hom_{stmod(\bG)}(C^\bu,D^\bu)$
if and only if the exact triangle $D^\bu \to cone(f-g)\to C[1]$ splits.   Using (\ref{eqn:Hom-Db-cx}), we 
conclude such a splitting exists if only if this exact triangle restricted to $\Gr$ admits a splitting for all $r > > 0$.
This, together with the isomorphism (\ref{eqn:Hom-Db-cx}) implies the injectivity of (\ref{eqn:bij-st}).  
The surjectivity of (\ref{eqn:bij-st}). follows from the
isomorphism (\ref{eqn:Hom-bij-cx}) and Theorem \ref{thm:St-rickard}.
\end{proof}

In the following proposition, 
  we do not rule out the existence of a chain of maps in $D^b(Mod(\bG))$
  beginning and ending with bounded complexes of finite dimensional
  modules $C^\bu \leftarrow X^\bu_1 \to X^\bu_2 \leftarrow \cdots \to D^\bu$
  which represents a map in $StMod(\bG)$ which is not a map in $stmod(\bG)$.

 \begin{prop}
 \label{prop:faithful}
 Let $\bG$ be a connected, irreducible linear algebraic group such that every non-trivial
 mock injective $\bG$ is infinite dimensional.  Then \ $stmod(\bG) \ \hookrightarrow  \ StMod(\bG)$ \ 
is a faithful embedding of tensor triangulated categories.
 \end{prop}
 
 \begin{proof} 
 We utilize the following commutative square for each $r > 0$
  \begin{equation}
 \label{eqn:compare-stable}
 \begin{xy}*!C\xybox{%
 \xymatrix{
stmod(\bG)) \ar[d] \ar[r] & StMod(\bG) \ar[d] \\
stmod(\Gr)) \ar[r] & StMod(\Gr) ). } \\
  }\end{xy}
 \end{equation}
 Observe that $stmod(\Gr)  \ \to \ StMod(\Gr)$ is fully faithful for all $r$ (utilizing the fact that
 maps in these categories are equivalence classes of maps of $\Gr$-modules).
 By Proposition \ref{prop:detect-faithful}, if $f \not=g \in Hom_{stmod(\bG))}(C^\bu,D^\bu)$,
 then $f, \ g$ have unequal images in $Hom_{stmod(\Gr)}(C^\bu_{|\Gr},D^\bu_{|\Gr})$ for $r > > 0$ and thus (by the
 previous observation) $f, \ g$ have unequal images in
 $Hom_{StMod(\Gr)}(C^\bu_{|\Gr},D^\bu_{|\Gr})$.   Consequently,  the commutativity of 
 (\ref{eqn:compare-stable}) implies that $stmod(\bG) \ \to \ StMod(\bG)$ \ 
is faithful.
  \end{proof}

 \vskip .1in
 
We summarize various categories we have considered, together with functors relating these categories,
in the following commutative diagram:
 \begin{equation}
 \label{eqn:multi-maps}
 \begin{xy}*!C\xybox{%
   \xymatrix{
\cI nj^b(Mod(\Gr)) \ar[r]  &    D^b(Mod(\Gr)) \ar[r]  & StMod(\Gr) \\
\cM ock^b(\bG) \ar[r]  \ar[u] &    D^b(Mod(\bG)) \ar[r]  \ar[u] &  StMod(\bG)  \ar[u] \\
mock^b(\bG) \ar[r]  \ar[u] \ar[d]&    D^b(mod(\bG)) \ar[r]  \ar[u] \ar[d] & stmod(\bG) \ar[d] \ar[u]\\  
 Perf(\Gr)) \ar[r] &    D^b(mod(\Gr)) \ar[r]  & stmod(\Gr).  }\\  
  }\end{xy}
 \end{equation}

    \vskip .1in
    
   The  functors established next with target $StMod(\bG)$ might be useful in future applications.

    \begin{prop}
    \label{prop:func}
    Let $f: (\bG, \cE) \to (\bG^\prime,\cE^\prime)$ be a map of algebraic groups of exponential type
    as in Definition \ref{defn:functor-exp}.
    Then restriction determines a well defined functor 
    $$f^*: StMod(\bG^\prime) \quad \to \quad  StMod(\bG)$$
    of tensor triangulated categories.  

    For any algebraic group $\bG$ of exponential type and any $r > 0$, induction determines a 
    well defined functor 
    $$ind_{\bG_{(r)}}^{\bG}: StMod(\bG_{(r)}) \quad \to \quad StMod(\bG)$$
    of triangulated categories. 
    
     If, in addition, $\bG$ is defined over $\bF_q$ for some $p$-th power $q$,
    then  induction determines a well defined functor 
    $$ind_{\bG(\bF_q)}^{\bG}: StMod(\bG(\bF_q)) \quad \to \quad StMod(\bG)$$
    of triangulated categories.
    \end{prop}
    
    \begin{proof}
    Since  $f^*: \cK^b(Mod(\bG^\prime)) \ \to \ \cK^b(Mod(\bG))$ is exact and preserves
    tensor products, it is a tensor triangulated functor.  By Proposition \ref{prop:further-cC},
    $f^*$ sends mock injective modules to mock injective modules, and thus
   determines  $f^*: StMod(\bG^\prime) \to StMod(\bG)$.
       
    We recall from \cite{CPS1} that a closed subgroup $H \subset \bG$ is called exact 
    if $ind_H^\bG: Mod(H) \to Mod(\bG)$ is exact (i.e., preserves short exact sequences).   
    By \cite[Thm4.3]{CPS1}, $\bG_{(r)} \hookrightarrow \bG$ is exact
    since $\bG/\Gr \ = \ \bG^{(r)}$ is affine.  Moreover, the induction functor $ind_{\bG_{(r)}}^{\bG}$
    sends (mock) injective $\Gr$-modules to injective $\bG$-modules.  
    We conclude that $ind_{\bG_{(r)}}^{\bG}$ induces the triangulated map of Verdier quotients 
    $$ind_{\bG_{(r)}}^{\bG}: StMod(\Gr) \simeq \cK^b(Mod(\Gr))/\cI nj^b(Mod(\Gr)) $$
    $$ \ \to \    \cK^b(Mod(\bG))/\cM ock^b(\bG) \simeq StMod(\bG).$$
    
This same proof justifies the tensor triangulated functor  $ind_{\bG(\bF_q)}^{\bG}: StMod(\bG(\bF_q)) \to StMod(\bG)$
upon replacing $\bG_{(r)} \hookrightarrow \bG$ by $\bG(\bF_q) \hookrightarrow \bG$.
    \end{proof}

    \vskip .2in


 \section{Support theory for complexes of $\bG$-modules}
 \label{sec:Gcomplexes}
 
 In this section, we extend our formulation of the support theory  $M \mapsto \Pi(\bG)_M$
 to bounded complexes $C^\bu$ of $\bG$-modules following our consideration of 
 various triangulated categories in the previous section.  As we show below, 
 \ $C^\bu \ \mapsto \ \Pi(\bG)_{C^\bu}$ \ depends only upon the isomorphism class of $C^\bu$
 in $StMod(\bG)$ and satisfies evident analogues of the good properties for $M \mapsto \Pi(\bG)_M$
 established in Theorem \ref{thm:properties}.
 This intertwining of support theory and stable module categories will be applied in 
 Section \ref{sec:stmod(bG)} to provide a classification of certain tensor triangulated
 subcategories of our stable module categories.
 
 We proceed by first formulating in Definition \ref{defn:Pi(G)-C} our support theory \ 
 $C^\bu \ \mapsto \ \Pi(G)_{C^\bu}$ for bounded complexes of $G$-modules
 with $G$ a finite group scheme.
Proposition \ref{prop:equalST} shows that this formulation corresponds to the usual
support theory for $G$-modules  under the homeomorphism of 
Theorem \ref{thm:St-rickard}.  After defining our support theory
 $C^\bu \ \mapsto \Pi(\bG)_{C^\bu}$ for bounded complexes of $\bG$-modules
 in Definition \ref{defn:Pi-Cdot}, we present in Theorem \ref{thm:stable-prop} 
 numerous good properties of this theory  whose proofs are derived from the 
 analogous properties of $M \mapsto \Pi(\Gr)_M$ using Proposition \ref{prop:equalST}.
 For $(\bG,\cE)$ an algebraic group of exponential type, we describe the subspace
 $\bP \cC(\bG)_{C^\bu} \ \subset \ \bP \cC(\bG)$ corresponding to $\Pi(\bG)_{C^\bu}
 \ \subset \Pi(\bG)$ in terms of ``actions at exponentials."

We recall the map of $k$-algebras $\epsilon_r: k\bG_{a(1)} \to k\bG_{a(r)}$ of (\ref{eqn:epsilon-r})
and observe that the restriction map $\epsilon_r^*: \cK^b(Mod(\bG_{a(r)}) \to \cK^b(Mod(\bG_{a(1)})$
is exact (i.e., a map of triangulated categories).  
Observe that the 0-module is a free module of rank 0 for any $G$ (for any affine group 
scheme $G$).  In particular, if $C^\bu \in CH^b(G)$ is contractible (i.e., homotopy equivalent to
the 0-complex), then $C^\bu$ is quasi-isomorphic to the 0-complex and thus to a complex 
of free $G$-modules.
     
  \vskip .1in
    
 \begin{defn}
\label{defn:Pi(G)-C}
Let $G$ be a finite group scheme over $k$ and let $C^\bu$ be a bounded complex of $G$-modules.
Let $\alpha_K: K[t]/t^p \to KG$ be a $\pi$-point whose equivalence  class $[\alpha_K]$ 
is a (Zariski) point of $\Pi(G)$.  Then $[\alpha_K]$ is said to be a point of the
$\Pi(G)_{C^\bu} \ \subset \ \Pi(G)$ if $\alpha_K^*(C^\bu_K)$ is not quasi-isomorphic 
to a bounded complex of free $K[t]/t^p$-modules.
 \end{defn}
 
 \vskip .1in

 The awkwardness of proving that $C^\bu \to \Pi(G)_{C^\bu}$ is well defined on objects
 of $StMod(G)$ is partially caused by the fact that $M \mapsto \Pi(G)_M$ is not 
 a well defined functor on $Mod(G)$.  For example, a map $M \to N$ of
 $G$--modules and a $\pi$-point $\alpha_K: K[t]/t^p \to KG$ might be such 
 that $\alpha_K^*(M)$ is not free (thus determining a point $\alpha_K$ of $\Pi(G)_M$) 
 whereas $\alpha_K^*(N)$ is free.
 
\begin{prop}
\label{prop:equalST}
Let $G$ be a finite group scheme.  
\begin{enumerate}
\item 
If a map $f:C^\bu \ \to  \ D^\bu$ of bounded complexes of $G$-modules satisfies the condition that
$cone(f)$  is quasi-isomorphic to a bounded complex of injective $G$-modules, 
then $\Pi(G)_{C^\bu}  \ = \ \Pi (G)_{D^\bu}$ (as subspaces of $\Pi(G)$).
\item
If $C^\bu,  \ D^\bu$ are bounded complexes of $G$-modules which are isomorphic in $D^b(Mod(G))/\cI nj^b(Mod(G))$,
then $\Pi(G)_{C^\bu}  \ = \ \Pi (G)_{D^\bu}$.
\item
For any $C^\bu$ in $D^b(Mod(G))$, the functor $\tilde \Phi: D^b(Mod(G)) \ \to \ StMod(G)$ of 
Theorem \ref{thm:St-rickard} satisfies the property that 
\begin{equation}
\label{eqn:Pi-Phi-C}
\Pi(G)_{C^\bu} \quad = \quad \Pi(G)_{\tilde \Phi(C^\bu)}.
\end{equation}
\item
For any $C^\bu$ in $D^b(Mod(G))$, $\Pi(G)_{C^\bu} = \emptyset$ if and only if $C^\bu$ is an object
of $\cI nj^b(Mod(G))$.
 \item
 If $C^\bu \to D^\bu \to E^\bu\to C^\bu[1] $ is an exact triangle in $D^b(Mod(G))$,
 then $\Pi(\bG)_{D^\bu} \ \subset \ \Pi(\bG)_{C^\bu} \cap \Pi(\bG)_{E^\bu}$.
\item
For any $C^\bu$ in $D^b(mod(G))$, $\Pi(G)_{C^\bu} \ \subset \ \Pi(\bG)$ is closed.
\end{enumerate} 
\end{prop}

\begin{proof}
Assertion (1) follows from Corollary \ref{cor:detect-iso}(3) since 
each $\pi$-point $\alpha_K: K[t]/t^p \to KG$ determines a triangulated map $\alpha_K^*:
D^b(Mod(G)) \to D^b(K[t]/t^p)$, preserving quasi-isomorphisms, commuting
with taking cones, and sending injective $G$-modules to free $K[T]/t^p$-modules.

Consider a ``roof" $s^{-1}\circ f$ representing a map from  $C^\bu$ to $D^\bu$ 
in $StMod(G)$; namely, a pair of maps in $D^b(Mod(G))$, \ 
$C^\bu \stackrel{s}{\leftarrow} E \stackrel{f}{\to} D^\bu$, \ with the cone of $s$ quasi-isomorphic 
to a bounded complex of injective modules.  Since maps in $StMod(G)$
are given by a calculus of fractions, we conclude that such a roof is an isomorphism in
$StMod(G)$ if and only $f$ is itself an isomorphism in $StMod(G)$.
Thus, assertion (2) follows from assertion (1).

Since the functor $\Psi$ of Theorem \ref{thm:rickard} is essentially surjective, it suffices 
to verify the identification (\ref{eqn:Pi-Phi-C}) for $C^\bu$ of the form $M[0]$ for some
$\bG$-module $M$.  This is immediate from the definitions.

Assertion (4) follows from assertion (3) and the ``projectivity test" for $G$-modules
which asserts that $\Pi(G)_M = \emptyset$ if and only if $M$ is a projective $G$-modules
\cite{P}.  Alternatively, (4) follows directly from Corollary \ref{cor:detect-iso}(1) and the 
fact that $\Phi$ is an equivalence of categories.

Assertion (5) follows from the exactness of $\Phi$ and the ``two out of three" property for
$M \mapsto \Pi(G)_M$ given in Theorem \ref{thm:properties}(4).

Finally, if $C^\bu$ in $D^b(mod(G))$, then $\tilde \Phi(C^\bu)$ is an object of $stmod(G)$
represented by a finite dimensional $G$-module.  Since $\Pi(G)_M\ \subset \ \Pi(G)$ is closed
for any finite dimensional $G$-module by Theorem \ref{thm:properties}(6), assertion (6) follows
 from assertion (3).
\end{proof}

The following definition of  $C^\bu \ \mapsto \ \Pi(\bG)_{C^\bu}$ is a natural extension 
of $M \mapsto \Pi(\bG)_M$ for a $\bG$-module $M$ as defined in Definition \ref{defn:Pi(bG)} and
 $C^\bu \ \mapsto \ \Pi(\Gr)_{C^\bu}$  
for a bounded complex of $\Gr$-modules as formulated in Definition  \ref{defn:Pi(G)-C}.

\begin{defn}
\label{defn:Pi-Cdot}
\vskip .2in
Let $\bG$ be a linear algebraic group and let $C^\bu$ be bounded complex of $\bG$-modules.
We define the subspace \ $\Pi(\bG)_{C^\bu} \ \subset \ \Pi(\bG)$ \ by 
$$\Pi(\bG)_{C^\bu} \quad \equiv \quad \varinjlim_r \Pi(\Gr)_{C^\bu_{|\Gr}} \quad \subset \quad 
\varinjlim_r \Pi(\Gr)\quad \equiv \quad \Pi(\bG)$$
where $\Pi(\Gr)_{C^\bu_{|\Gr}} \subset \ \Pi(\Gr)$ is defined in Definition \ref{defn:Pi(G)-C}.
\end{defn}

 \vskip .1in
 
 We  state and prove our generalization of Theorem \ref{thm:properties} to bounded 
 complexes of $\bG$-modules, thereby establishing the basic properties of the
 support theory $C^\bu \to \Pi(\bG)_{C^\bu}$.

\begin{thm}
 \label{thm:stable-prop}
Let $\bG$ be a linear algebraic group over $k$.  Then 
 $$C^\bu \in CH^b(Mod(\bG)) \quad \ \text{maps  to} \quad \Pi(\bG)_{C^\bu} \ \subset \ \Pi(\bG)$$
as defined in Definition \ref{defn:Pi-Cdot} satisfies the following properties.
\begin{enumerate}
\item
$\Pi(\bG)_{C^\bu} \ = \ \emptyset$ \ if and only \ $C^\bu \ \in \ \cM ock^b(\bG)$.
 \item
 $\Pi(\bG)_{C^\bu}  \quad = \quad \Pi(\bG)_{C^\bu[n]}$ as\ subsets of $V(\bG)$
 \item
 If $C^\bu, \ D^\bu$ are bounded complexes of $\bG$-modules which are isomorphic
 in $StMod(\bG)$, then $\Pi(\bG)_{C^\bu} \ = \ \Pi(\bG)_{D^\bu}$.
  \item
 If $C^\bu \to D^\bu \to E^\bu\to C^\bu[1] $ is an exact triangle in $StMod(\bG)$,
 then $\Pi(\bG)_{D^\bu} \ \subset \ \Pi(\bG)_{C^\bu} \cap \Pi(\bG)_{E^\bu}$.
\item
If $C^\bu \ = \ \bigoplus_{i \in I} C^\bu_i$, then \
    $\Pi(\bG)_{C^\bu} \ =  \ \bigcup_{i\in I} \Pi(\bG)_{C^\bu_i}$.
\item
 $\Pi(\bG)_{C^\bu \otimes D^\bu} \ = \ \Pi(\bG)_{C^\bu} \cap \Pi(\bG)_{D^\bu}$.
\item
If $C^\bu \in CH^b(mod(\bG))$, then \ $\Pi(\bG)_{C^\bu}\ \subset \ \Pi(\bG)$ \ is closed.
\end{enumerate}
\end{thm} 
    
\begin{proof}
 As defined in Definition \ref{defn:cMock}, $C^\bu$ is an object of $Mock^b(\bG)$ 
 if and only if its restriction to each $\Gr$ lies in $Inj^b(Mod(\Gr)$) which is equivalent by 
 Proposition \ref{prop:equalST}(4) to the condition that each $\Pi(\Gr)_{C^\bu}$ is empty. This proves the first assertion

 The second assertion follows from the fact that $(\cE_{\Lambda_r(\ul B)}\circ \epsilon_r)^*(-)$
 commutes with shifts $(-)[n]$ of chain complexes and the fact that  $(-)[n]$ preserves quasi-isomorphisms.

Proposition \ref{prop:equalST}(2) implies assertion 3.
Assertion (4) follows from the exactness of restriction of $\bG$-modules to $\Gr$-modules
and Proposition \ref{prop:equalST}(5).  
Similarly, assertions (5) and (6) follow from the exactness of restriction of $\bG$-modules 
to $\Gr$-modules,  Proposition \ref{prop:equalST}(3), and the corresponding properties for $M \to \Pi(\bG)_M$
given in Theorem \ref{thm:properties}.

Since $\Pi(\bG)$ is equipped with the colimit topology, $\Pi(\bG)_{C^\bu}\ \subset \ \Pi(\bG)$ is closed
if and only if each $\Pi(\Gr)_{C^\bu}\ \subset \ \Pi(\Gr)$ is closed.  Thus, assertion (7) follows from
Proposition \ref{prop:equalST}(6).
\end{proof}

\vskip .1in

\begin{remark}
\label{rem:real}
It is natural to ask what if any new closed subsets of $\Pi(\bG)$ (other than those of the
form $\Pi(\bG)_{M[0]}$ for $M$ finite dimensional) are realized as $\Pi(\bG)_{C^\bu}$ 
for $C^\bu$ a bounded complex of finite dimensional $\bG$-modules.  

A serious impediment to understanding what subsets are realized as $\Pi(\bG)_{C^\bu}$
is that the functor $\Phi_{(r)}$ for a given $r > 0$ (as in Theorem \ref{thm:St-rickard}) depends
upon the triangulated structure of $StMod(\Gr)$ which is not natural with respect to passing
from $\Gr$ to $\bG_{(r+1)}$.
Thus, if we take a bounded complex of finite dimensional $\bG$-modules $C^\bu$, we 
may associate for each $r > 0$ a finite dimensional $\Gr$-module $\Phi_{(r)}(C^\bu)$,
but we lack a method to construct a (finite dimensional) $\bG$-module using the
family $\{ \Phi_{(r)}(C^\bu) , \ r >0 \}$.
\end{remark}

We introduce the following terminology for ``realizable" subsets of $\Pi(\bG)$.

\begin{terminology}
\label{term:realizable}
Let $\bG$ be a  linear algebraic group and consider a subspace $X \ \subset\ \Pi(\bG)$. 
\begin{itemize}
\item
$X$ is said to be $Mod(\bG)$-realizable (respectively, $mod(\bG)$-realizable) if there exists 
some $\bG$-module (resp., finite dimensional $\bG$-module) $M$ such that $X \ = \ \Pi(\bG)_M$.
\item
$X$ is said to be $StMod(\bG)$-realizable (respectively, $stmod(\bG)$-realizable) if there exists 
some bounded complex of $\bG$-modules (resp., bounded complex of finite dimensional $\bG$-modules) $C^\bu$ 
such that $X \ = \ \Pi(\bG)_{C^\bu}$.
\item
$X$ is said to be locally $StMod(\bG)$-realizable (respectively, locally $stmod(\bG)$-realizable) if there exists 
some collection of bounded complexes of $\bG$-modules (resp., collection of bounded complexes
of finite dimensional $\bG$-modules) $\{C^\bu_{\alpha} \}$ 
such that $X \ = \ \bigcup_\alpha \Pi(\bG)_{C^\bu_\alpha}$.
\end{itemize}
\end{terminology}

\vskip .1in

	The identification (\ref{eqn:cC-iso})  of $M \to \Pi(\bG)_M$ with  $M \to \bP\cC (\bG)_M$ for $(\bG,\cE)$ an
algebraic group of exponential type provides a more ``concrete" interpretation of $M \to \Pi(\bG)_M$.
We proceed to show that if  $(\bG,\cE)$ an algebraic group of exponential type, then this identification extends to
bounded complexes of $\bG$-modules.   

\begin{defn}
\label{defn:PC(G)-C}
Consider an algebraic group $(\bG,\cE)$ of exponential type.  For any bounded complex $C^\bu$ of $\bG$-modules,
we define the subspace \ $\bP\cC(\bG)_{C^\bu} \quad \subset \quad \bP\cC(G)$ by identifying its $K$-points for any
field extension $K/k$ to consist of those $K$-points of $\bP\cC(G)$ represented by a $K$-point $\ul B \in \cC(\bG)$
such that 

 $(\cE_{\Lambda_r(\ul B)} \circ \epsilon_r)^*(C^\bu)$ \ is not quasi-isomorphic to a complex of free modules \\
for any $r$ sufficiently large that  $B_i = 0$ for $i \geq r$.

We define $\bP \cC_r(\bG)_{C^\bu}$ as the intersection \ $\bP \cC_r(\bG) \cap \bP\cC(\bG)_{C^\bu}$, so that 
$$\bP\cC(\bG)_{C^\bu} \quad = \quad \bigcup_r \bP \cC_r(\bG)_{C^\bu}.$$
\end{defn}

\vskip .1in

\begin{prop}
 \label{prop:expC}
 Let $(\bG, \cE)$ be an algebraic group of exponential type, fix $r > 0$, and let $C^\bu$ be a
a bounded complex of $\bG$-modules.  The homeomorphism \\
 $\Phi \circ \Psi \circ (\rho_r \circ \lambda): \bP \cC_r(\bG) \ \to \ \Pi(\Gr)$ of Theorem \ref{thm:cC-Pi} restricts
 to a homeomorphism
 \begin{equation}
 \label{eqn:cC-iso-C}
 \Phi \circ \Psi \circ (\rho_r \circ \lambda): \bP \cC_r(\bG)_{C^\bu} \ \stackrel{\sim}{\to} \ \Pi(\Gr)_{C^\bu_{|\Gr}},
 \end{equation}
 thereby determining the homeomorphism \ $\bP \cC(\bG)_{C^\bu} \ \stackrel{\sim}{\to} \ \Pi(\bG)_{C^\bu}.$
\end{prop}

\begin{proof}
To prove (\ref{eqn:cC-iso-C}), we must compare Definition \ref{defn:Pi(G)-C} with Definition \ref{defn:PC(G)-C}.
Exactly as in the proof of Theorem \ref{thm:cC-Pi}, this comparison is made by
juxtaposing the determination of $\Psi \circ (\rho_r \circ \lambda_r): \bP \cC_r(\bG) \to \bP |\Gr|$
in Theorem \ref{thm:rho-lambda}(4), the
definition of $\rho_r$ in Definition \ref{defn:rho}, and the determination of $\phi \circ \Psi: \bP V_r(\Gr) \to \Pi(\Gr)$
in Theorem \ref{thm:1par,pi}.
\end{proof}

The following proposition is the extension to complexes of Proposition \ref{prop:exp-deg}(1).  The proof
follows immediately upon recalling that the $\pi$-point associated to the 1-parameter subgroup $\cE_{\ul B}$ with
$B_i = 0, \ i \geq r$ is $\alpha_K \equiv \cE_{ \Lambda_r(\ul B)} \circ \epsilon_r: K[t]/t^p \to K\bG$ so that the
conditions of Definition \ref{defn:Pi(G)-C} and Definition \ref{defn:PC(G)-C} match in light of the identification $\Pi(\bG)_{C^\bu}
\ \equiv \ \varinjlim_r \Pi(\Gr)_{C^\bu}$ of Definition \ref{defn:Pi-Cdot}.

\begin{prop}
Let $(\bG, \cE)$ be an algebraic group of exponential type and let $C^\bu$ be a
a bounded complex of $\bG$-modules.
If $C^\bu$ is quasi-isomorphic to a bounded complex each of whose terms has 
bounded exponential degree (as in Proposition \ref{prop:exp-deg}, then 
$$\bP\cC(\bG)_{C^\bu} \ = \ \pi_r^{-1}(\bP\cC_r(\bG)_{C^\bu}), \quad r > > 0.$$
\end{prop}

In anticipation of the classification of certain subcategories in Section \ref{sec:stmod(bG)},
we remark that the full subcategory of $stmod(\bG)$ (respectively.  $StMod(\bG)$) consisting
of bounded complexes of $\bG$-modules whose terms have bounded exponential degree is a
thick triangulated subcategory (resp., localizing subcategory), but not a tensor ideal.
\vskip .2in


 \section{Classifying subcategories of $stmod(\bG)$ and $StMod(\bG)$}
 \label{sec:stmod(bG)}
  
 In this section, we use the classification of thick, tensor ideals of $stmod(\Gr)$
 (given, for example, in \cite{R1}) to classify the {\it $(r)$-complete} thick tensor ideals
 of $stmod(\bG)$ in Theorem \ref{thm:HPS}.   The formulation of the property
`` $(r)$-complete" for a thick tensor ideal $\cC \subset stmod(\bG)$ is introduced in Definition 
\ref{defn:Cr}.  This property is naturally suggested by the relationship of
the support varieties for $\bG$-modules and those associated to restrictions of
$\bG$-modules to Frobenius kernels $\{ \Gr, \ r > 0 \}$.  We also present an
analogous classification for the {\it $(r)$-complete localizing} tensor 
triangulated subcategories of $StMod(\bG)$.

We remind the reader that a full triangulated subcategory $\cC$ of a triangulated
 category $\cD$ is a {\it thick subcategory} if and only if every object of $\cD$ which is a direct summand of 
 an object of $\cC$ is itself an object of $\cC$ \cite{R1}.    If the triangulated category 
 $\cD$ has a (for convenience, symmetric)
 tensor structure, then  a full triangulated subcategory $\cC$
 of $\cD$ is said to be a {\it tensor ideal} if tensoring any object of $\cC$ with an object of $\cD$ is again
 an object of $\cC$.  The thick
tensor ideal in $\cC$ generated by a collection of objects of $\cC$ is the full triangulated subcategory of $\cC$
whose objects are obtained by repeatedly applying the operations of taking finite sums of objects, taking summands of 
objects, taking cones of maps, and taking tensor products with arbitrary objects of $\cC$.

 \begin{defn}
 \label{defn:Cr}
 Consider a linear algebraic group $\bG$ and let $\cC \subset stmod(\bG)$ be a 
 triangulated subcategory.  
 For each $r > 0$, we denote by $\cC_{|\Gr} \ \subset \ stmod(\Gr)$ the essential image of  $\cC$
 under the restriction functor $stmod(\bG) \to stmod(\Gr)$.  In other words, \ $\cC_{|\Gr}$ \ is the 
 full subcategory of $\cC$ whose objects are those objects of $\cC$ isomorphic to objects obtained
 by restriction of objects in $stmod(\bG)$.
 
 If  $\cC  \ \subset \ stmod(\bG)$ is a thick tensor ideal, then we denote by  $\cC_{(r)}$  the 
thick tensor ideal of $stmod(\Gr)$ generated by $\cC_{|\Gr}$.

We say that a thick tensor ideal $\cC \subset stmod(\bG)$ is {\it $(r)$-complete} if the following 
condition is satisfied for every $C^\bu \in stmod(\bG)$: 
 \begin{equation}
 \label{eqn:saturate}
C^\bu \in \cC \quad \iff \quad C^\bu_{|\Gr} \in \ \cC_{(r)}, \  \forall r > 0.
\end{equation}
 \end{defn}
 
 \vskip .1in
 We remark that even if $\cC$ is $r$-complete, $\cC_{|\Gr}$ is unlikely to be a tensor
ideal in $stmod(\Gr)$ since typically there are objects of $stmod(\Gr)$ which are not 
restrictions of objects of $stmod(\bG)$. 
 
 \begin{defn}
 \label{defn:cat-space}
 Consider a linear algebraic group $\bG$.
  For any collection $\cS$ of bounded complexes of finite dimensional $\bG$-modules, we define 
 the locally closed subset 
    $$\Pi(\bG,\cC) \ \equiv \ \bigcup_{C^\bu \in \cS} \Pi(\bG)_{C^\bu} \  \subset \ \Pi(\bG)$$
and define 
    $$\Pi(\Gr,\cC_{| \Gr})  \ \equiv \ \bigcup_{C^\bu \in \cS} \Pi(\Gr)_{C^\bu_{| \Gr}} \  \subset \ \Pi(\Gr),$$
    where $C^\bu_{| \Gr}$ is the restriction of $C^\bu$ to $\Gr$.
    
    For any subset  $X  \subset  \Pi(\bG)$, we define the full subcategory
    $$\cC_X  \ \equiv \ \langle \{ C^\bu \in stmod(\bG), \ such \ that \ \bP V(\bG)_{C^\bu} \subset X\} \rangle \ \subset \ stmod(\bG).$$ 
    \end{defn}

  \vskip .1in
   
 \begin{prop}
 \label{prop:saturation}
  Let $\bG$ be a linear algebraic group and let 
 $\cC \ \subset \ stmod(\bG)$ be a thick tensor ideal.  We define  the $(r)$-completion of $\cC$,
 \ $\cC^\vee$, \ to be the full subcategory of $stmod(\bG)$ whose objects are
 those bounded complexes  $C^\bu$ of finite dimensional $\bG$-modules such
 that  the restriction of $C^\bu$ to $\Gr$ lies in $\cC_{(r)}$ for every $r > 0$.
 
 Then $\cC \ \mapsto \ \cC^\vee$ satisfies the following properties.
 \begin{enumerate}
\item
$\cC^\vee$ is a thick tensor ideal of  $stmod(\bG)$.
\item
$\cC_{(r)} \ = \ (\cC^\vee)_{(r)}$ for all $r > 0$, so that $(\cC^\vee)^\vee \ = \ \cC^\vee$.
\item
$\cC^\vee$ is the the minimal $(r)$-complete thick tensor ideal 
of $stmod(\bG)$ containing $\cC$.
\item 
The natural embedding 
$\Pi(\Gr,\cC_{| \Gr}) \ \hookrightarrow \Pi(\Gr,\cC_{(r)})$
is the identity for each $r > 0.$

\end{enumerate}
 \end{prop}
 
 \begin{proof}
 If $C^\bu \to D^\bu \to E^\bu \to C^\bu[1]$ is an exact triangle in $stmod(\bG)$ with 
 $C^\bu \to D^\bu$ in $\cC^\vee$, then the restriction of this exact triangle to
 $stmod(\Gr)$ is an exact triangle in $\cC_{(r)}$ for all $r > 0$. 
 This implies that $E^\bu$ is an object of $\cC^\vee$,  Similarly, if $C^\bu \in C^\vee$
 and $X^\bu$ is an arbitrary bounded complex of finite dimensional $\bG$-modules,
 then any summand of $C^\bu$ is an object of $C^\vee$ and $C^\bu \otimes X^\bu$ is
 an object of $C^\vee$ for any $X^\bu \in stmod(\bG)$,
 since each $\cC_{(r)} \subset stmod(\Gr)$ is a thick tensor ideal.  Thus,
 $\cC^\vee$ is a thick tensor ideal of $stmod(\bG)$.
 
The operations on objects of $\cC$ which produce objects of $\cC^\vee$ restrict
to internal operations on $\cC_{(r)}$, so that $\cC_{(r)} \ = \ (\cC^\vee)_{(r)}$.
This implies that $(\cC^\vee)^\vee \ = \ \cC^\vee$ and that 
$\cC^\vee$ is $(r)$-complete.

Observe that objects of $\cC_{(r)}$ are obtained by starting with objects of 
$\cC_{| \Gr} \ \subset \ stmod(\Gr)$ and successively applying
the operations of taking cones, taking direct summands, and tensoring with objects of $stmod(\bG)$.  
By Theorem \ref{thm:stable-prop}, these operations 
preserve the property that if the support of the input object of an operation is contained in 
$\Pi(\Gr,\cC_{|\Gr})$, then the support of the output object is also contained in $\Pi(\Gr,\cC_{|\Gr})$.
Consequently,  \ $\Pi(\Gr,\cC_{| \Gr}) \ = \ \Pi(\Gr,\cC_{(r)})$.
 \end{proof}
    
     \vskip .1in
 
The definition of an $(r)$-complete thick tensor ideal of $stmod(\bG)$ was
made in anticipation of the following result.
 
 \begin{prop} 
 \label{prop:sat-cX}
  Let $\bG$ be a linear algebraic group  and let $X \ \subset \ \Pi(\bG)$ 
be a subspace.  Then $\cC_X$ is a thick tensor ideal of $stmod(\bG)$ which  
is $(r)$-complete.
 \end{prop}
 
 \begin{proof}
 Theorem \ref{thm:stable-prop}(4) easily implies that $\cC_X$ is a thick triangulated subcategory
 of $stmod(\bG)$; consequently, Theorem \ref{thm:stable-prop}(6) implies that $\cC_X$ is a 
 thick tensor ideal in $stmod(\bG)$.
 
 Assume that $C^\bu$ is an object of $\cC_X^\vee$; in other words, $C^\bu$
 is a bounded complex of finite dimensional $\bG$-modules whose restriction
 to each $\Gr$ lies in $(\cC_X)_{(r)}$.  By Proposition \ref{prop:saturation}(4) and the evident inclusion
 $\Pi(\Gr,(\cC_X)_{| \Gr}) \ \subset \ X$, we conclude that $\Pi(\Gr,(\cC_X)_{(r)}) \ \subset \ X$ for all $r > 0$,
 so that $\Pi(\bG)_{C^\bu} \ \subset \ X$ and thus $C^\bu$ is an object of $ \cC_X$.
\end{proof}

  We proceed to show that  locally $stmod(\bG)$-realizable subsets of $\Pi(\bG)$ as in Terminology \ref{term:realizable}
classify $(r)$-complete, thick tensor ideals of  $stmod(\bG)$.  Our proof is heavily dependent upon
 the classification of thick tensor ideals of
 $stmod(\Gr)$ first established in \cite{BCR} for finite groups and then for general finite group schemes in \cite{FP2}.
  
    \begin{thm}
    \label{thm:HPS}
    Let $\bG$ be a linear algebraic group.   The correspondences 
    $$X \ \mapsto \ \cC_X \quad \cC \ \mapsto \ \Pi(\bG,\cC).$$
restrict to give mutually inverse bijections
    $$\{  locally \ stmod(\bG){\text-}realizable\ subsets \ X \subset \Pi(\bG) \}  \ \longleftrightarrow  $$
   $$ \{(r){\text-}complete, \ thick \ tensor  \ ideals \  \cC \subset stmod(\bG) \}.$$
    \end{thm}
    
    \begin{proof}
    Observe that  $\Pi(\bG,\cC_X) \ \subset \ X$ for any subspace $X \subset \Pi(\bG)$.  Assume now 
    that $X$ is  a locally $stmod(\bG)$ realizable subspace of $\Pi(\bG)$.  Then for any point 
    $x \in X$, there is there exists some bounded complex  $C^\bu_x$ of finite dimensional $\bG$-modules
    with $\Pi(\bG)_{C^\bu_x} \ \subset \ X$ such that $x \in \Pi(\bG)_{C^\bu_x}$.  Hence, $X \subset \ \Pi(\bG,\cC_X)$ 
    and thus $\Pi(\bG,\cC_X) \ = \ X$
    
  To complete the proof, we prove that the evident inclusion $\cC \ \subset \ \cC_{\Pi(\bG,\cC)}$ 
 of full, triangulated subcategories of $stmod(\bG)$ is an equivalence if $\cC \subset \stmod(\bG)$ is 
an $(r)$-complete tensor ideal.   Assume that $\cC$ is a thick tensor ideal and consider some object $E^\bu$ of $\cC_{\Pi(\bG,\cC)}$.   
 Since $\cC_{|\Gr} \ \subset \ \cC_{(r)}$, the restriction of $E^\bu$ lies in $\cC_{\bP V(\Gr,\cC_{(r)})}$.   
We apply the classification of thick tensor ideals of $stmod(\Gr)$, given in \cite[Thm6.3]{FP2} 
based upon the construction of Rickard idempotents in \cite{R2}, which tells us for the thick tensor ideal
$\cC_{(r)} \ \subset \ stmod(\Gr)$ that 
$\cC_{\bP V(\Gr,\cC_{(r)})} = \cC_{(r)}$ as subcategories of $stmod(\Gr)$.  Thus, the restriction of $E^\bu$ to 
each $\Gr$ lies in $\cC_{(r)}$ so that $E^\bu $ is an object of $\cC^\vee$.  Hence, if $\cC$ is also $(r)$-complete,
then $E^\bu \in \cC$.
\end{proof}
     
 \vskip .1in

    Rickard constructs idempotent endofunctors on $StMod(\Gr)$,
  \  $ \cE_{\cC}(-)$ and  $\cF_{\cC}(-)$, associated to a thick subcategory $\cC \subset stmod(\Gr)$.
  Among other properties, these functors satisfy the condition
 for any $X \in StMod(\Gr)$ that there is a natural exact triangle
 $$E_{\cC_{(r)}}(X) \to X \to F_{\cC_{(r)}}(X) \to E_{\cC_{(r)}}(X)[1]$$
 isomorphic to 
 $$E_{\cC_{(r)}}(k) \otimes X \to X \to F_{\cC_{(r)}}(k)\otimes X \to (E_{\cC_{(r)}}(k)\otimes X)[1].$$
   If $\cC$ is a tensor ideal inside $stmod(\Gr)$, then Rickard 's arguments using  
    adjunction relating $Hom(M^\vee,\cF_{\cC}(k))$ to  $M \otimes \cF_{\cC}(k)$ (and similarly
    for $\cE_{\cC}(k))$) provide a test for whether or not $M \in stmod(\Gr)$ 
    belongs to $\cC$.  This argument requires the compactness of $M$ 
    and the tensor ideal condition on $\cC$.

     We complement Theorem \ref{thm:HPS} with the following extension of Rickard's test
in terms of Rickard's idempotents, applying now to $C^\bu$ an object of $stmod(\bG)$ and 
$\cC \subset stmod(\bG)$ an $(r)$-complete, thick tensor ideal.

 \begin{thm}
 \label{thm:belong}
 Let $\bG$ be a  linear algebraic group and let $\cC \ \subset \ stmod(\bG)$ be an 
thick tensor ideal which is $(r)$-complete.  For each $r > 0$, we consider the Rickard idempotent functors 
 $$E_{\cC_{(r)}}(-), \  F_{\cC_{(r)}}(-): StMod(\Gr) \to StMod(\Gr)$$   
 for the  thick tensor ideal $\cC_{(r)} \subset stmod(\Gr)$ of Definition \ref{defn:Cr}
 associated to $\cC$.
  \begin{enumerate}
  \item
   $C^\bu \ \in \ \cC$ if and only if for all $r > 0$
\begin{equation}
\label{eqn:F}
F_{\cC_{(r)}}(C^\bu_{|\Gr}) \ \simeq \ F_{\cC_{(r)}}(k) \otimes C^\bu_{|\Gr} \  = 0 \quad {\text in} \ StMod(\Gr).
\end{equation}
  \item
  $C^\bu \ \in \ \cC$ if and only if for all $r > 0$
 \begin{equation}
\label{eqn:E}
E_{\cC_{(r)}}(C^\bu_{|\Gr}) \ \simeq \ E_{\cC_{(r)}}(k) \otimes C^\bu_{|\Gr} \  \simeq \ C^\bu_{|\Gr}
 \quad {\text in} \  StMod(\Gr).
\end{equation}
 \end{enumerate}   
\end{thm}

\begin{proof}
We assume that $\cC \ \subset \ stmod(\bG)$ is an $(r)$-complete thick tensor ideal.
If $C^\bu \in \cC$, then $C^\bu_{|\Gr} \in \cC_{(r)}$ for all $r > 0$ so that by Rickard's results
$F_{\cC_{(r)}}(C^\bu_{|\Gr}) = 0 \in StMod(\Gr)$ for all $r > 0$ (see \cite{FP2}).  Conversely,
if $F_{\cC_{(r)}}(C^\bu_{|\Gr}) = 0 \in StMod(\Gr)$ for all $r > 0$, then Rickard's results for $\Gr$ 
tell us that  $C^\bu_{|\Gr} \in \cC_{(r)}$ for all $r > 0$; since $\cC \ \subset \ stmod(\bG)$ is 
$(r)$-complete, this implies that $C^\bu \in \cC$.

 If $C^\bu \in \cC$, then Rickard's results tell us
 that $E_{\cC_{(r)}}(C^\bu_{|\Gr}) \  \simeq \  C^\bu_{|\Gr} \in StMod(\Gr)$ for all $r > 0$.  Conversely, 
 if $E_{\cC_{(r)}}(C^\bu_{|\Gr}) \  \simeq\  C^\bu_{|\Gr}$  in $StMod(\Gr)$ for all $r > 0$, then Rickard's
 exact triangles
$$E_{\cC_{(r)}}(C^\bu_{|\Gr}) \to C^\bu_{|\Gr} \to F_{\cC_{(r)}}(C^\bu_{|\Gr}) \to E_{\cC_{(r)}}(C^\bu_{|\Gr})[1]$$
imply that $F_{\cC_{(r)}}(C^\bu_{|\Gr}) = 0 \in StMod(\Gr)$ for all $r > 0$
so that (\ref{eqn:F}) implies that $C^\bu \in \cC$.
\end{proof}
    
    \vskip .1in
    
       We next make explicit how our earlier definitions and constructions for subcategories $\cC \ \subset stmod(\bG)$
can be modified to yield similar definitions and constructions for subcategories $\tilde \cC \ \subset StMod(\bG)$.
We remind the reader that a  {\it localizing subcategory} of a triangulated category admitting arbitrary direct sums
is a full triangulated subcategory closed under isomorphisms and arbitrary direct sums.
    
 In parallel with Definition \ref{defn:Cr}, we denote by  $\tilde \cC_{|\Gr}$ the essential image 
 of a localizing subcategory $\tilde C \ \subset \ StMod(\bG)$ under the restriction map 
 $StMod(\bG) \to StMod(\Gr)$; we denote by 
$(\tilde \cC)_{(r)}^\oplus\ \subset \ StMod(\Gr)$ the localizing subcategory generated by  $\tilde \cC_{|\Gr}$.
 We say that a localizing subcategory 
$\tilde C \ \subset \ StMod(\bG)$ is an  {\b\it$(r)$-complete localizing}  subcategory of $StMod(\bG)$ if
\begin{equation}
 \label{eqn:saturate-Mod}
C^\bu \in \tilde\cC \quad \iff \quad C^\bu_{|\Gr}  \  \in \ (\tilde \cC)_{(r)}^\oplus, \  \forall r > 0.
\end{equation}
As for subcategories of $stmod(\bG)$, we define $\Pi(\bG,\tilde \cC)\ \subset \ \Pi(\bG)$ for any subcategory 
$\tilde C \ \subset \ StMod(\bG)$ to be the subset 
    $$\Pi(\bG,\tilde\cC) \ \equiv \ \bigcup_{C^\bu \in \tilde\cC} \Pi(\bG)_{C^\bu} \  \subset \ \Pi(\bG).$$

We provide a natural analogue of Proposition \ref{prop:saturation}, with $(r)$-complete localizing
subcategories of $StMod(\bG)$ playing the role $(r)$-complete thick ideals of $stmod(\bG)$.

\begin{prop}
 \label{prop:saturation-Mod}
  Let $\bG$ be a linear algebraic group and let 
 $\tilde \cC \ \subset \ StMod(\bG)$ be a triangulated subcategory.  We denote by $\tilde \cC^\oplus$ 
 the full subcategory of $StMod(\bG)$ whose objects are
 those bounded complexes  $C^\bu \in Stmod(\bG)$ with the property
 that  the restriction of $C^\bu$ to $\Gr$ lies in $(\tilde \cC)_{(r)}^\oplus$ for every $r > 0$.
 Then $\tilde \cC^\oplus \ \subset StMod(\bG)$ satisfies the following properties:
 \begin{enumerate}
 \item
$\tilde \cC^\oplus$ is a localizing subcategory of $StMod(\bG)$. 
\item
$(\tilde \cC)_{(r)}^\oplus \ = \ (\tilde \cC^\oplus)_{(r)}^\oplus$, 
so that $(\tilde \cC^\oplus)^\oplus \ = \ \tilde \cC^\oplus$.
\item
$\tilde \cC^\oplus$ is the minimal $(r)$-complete localizing subcategory of $StMod(\bG)$ containing $\tilde \cC$.
\item
The natural embedding $\Pi(\Gr,\tilde \cC_{|\Gr}) \ \hookrightarrow \ \Pi\Gr,\tilde \cC^\oplus)$ is the identity
for each $r > 0$.
\end{enumerate}
 \end{prop}
 
 \begin{proof}
 The proof is a repetition of the proof of Proposition \ref{prop:saturation},
replacing $stmod(\bG)$ by $StMod(\bG)$, replacing $(r)$-complete thick ideals by $(r)$-complete localizing
subcategories, and replacing $(-)^\vee$ by $(-)^{\oplus}$.  

 In more detail, if $C^\bu \to D^\bu \to E^\bu \to C^\bu[1]$ is an exact triangle in $StMod(\bG)$ with 
 $C^\bu \to D^\bu$ in $\tilde \cC^\oplus$, then the restriction of this exact triangle to
 $StMod(\Gr)$ is an exact triangle in $\tilde \cC_{(r)}^\oplus$ for all $r > 0$ so that $E^\bu$ is in $\tilde \cC^\oplus$.   
 Similarly, if  $\{ C^\bu_\alpha, \ \alpha \in A\}$
 is a family of bounded complexes in $\tilde \cC^\oplus \subset StMod(\bG)$, then the restriction 
 of $\bigoplus_{\alpha \in A} C^\bu_\alpha$ to $\Gr$ lies in the localizing subcatory $(\tilde \cC)_{(r)}^\oplus \
 \subset \ StMod(\Gr)$ for every $r > 0$ so that $\bigoplus_{\alpha \in A} C^\bu_\alpha$
 is also in $\tilde \cC^\oplus$.  Thus, $\tilde \cC^\oplus \subset StMod(\bG)$ is a localizing subcategory.
 
The operations on objects of $\tilde \cC$ (taking cones and arbitrary direct sums) which produce 
objects of $\tilde \cC^\oplus$ beginning with objects of $\tilde \cC$ restrict
to internal operations on $(\tilde \cC)^\oplus_{(r)}$ since this is a localizing subcategory, 
so that $(\tilde \cC)_{(r)}^\oplus \ = \ (\tilde \cC^\oplus)_{(r)}^\oplus$.
This immediately implies that $(\tilde \cC^\oplus)^\oplus \ = \ \tilde \cC^\oplus$ and that 
$\tilde \cC^\oplus$ is an $(r)$-complete localizing subcategory of $StMod(\bG)$.

Finally, observe that  the operations of taking cones and arbitrary direct sums in $StMod(\Gr)$ 
preserve the property that if the support of the input object of an operation is contained in 
$\Pi(\Gr,\tilde\cC_{|\Gr})$, then the support of the output object is also contained in $\Pi(\Gr,\tilde\cC_{|\Gr})$.
As explained for proof of Proposition \ref{prop:saturation}(4), this implies the last assertion.
 \end{proof}
 
The following definition is the evident analogue for subcategories of $StMod(\bG)$ of 
Definition \ref{defn:cat-space}.  

 \begin{defn}
    \label{defn:subsets-arbitrary}
 Consider a linear algebraic group $\bG$.  For any subcategory $\tilde \cC \ \subset \ StMod(\bG)$,
 we define
$$ \Pi(\Gr,\tilde \cC_{| \Gr})  \ \equiv \ \bigcup_{C^\bu \in \tilde\cC} \Pi(\Gr)_{C^\bu_{| \Gr}} \quad
    \Pi(\bG,\tilde \cC) \ \equiv \ \bigcup_{C^\bu \in \tilde \cC} \Pi(\bG)_{C^\bu}.$$
 
 For any  subset $X$ of  $\Pi(\bG)$, we define the full subcategory
    $$\tilde \cC_X  \ \equiv \ \langle \{ C^\bu \in StMod(\bG) \ such \ that \ \Pi(\bG)_{C^\bu} \subset X\} \rangle \ \subset \ StMod(\bG).$$
\end{defn}
    
\vskip .1in
    
 We have the following analogue of Proposition \ref{prop:sat-cX} for $\tilde \cC_X \subset StMod(\bG)$.
 
 \begin{prop}
 \label{prop:sat-local}
 Let $\bG$ be a  linear algebraic group and let  $X \subset \Pi(\bG)$ be a subset.
 Then $\tilde \cC_X$ is an $(r)$-complete localizing tensor ideal of $StMod(\bG)$.  
 \end{prop}
 
 \begin{proof}
 By Theorem \ref{thm:stable-prop}(4),(5), \ $\tilde \cC_X$ is a localizing tensor ideal of $StMod(\bG)$. 
 
  Assume that $C^\bu$ is an object of $\tilde \cC_X^\oplus$; in other words, $C^\bu$
 is a bounded complex of $\bG$-modules whose restriction
 to each $\Gr$ lies in $(\tilde C_X)_{(r)}^\oplus$.  By Proposition \ref{prop:saturation}(4) and the evident inclusion
 $\Pi(\Gr,(\tilde \cC_X)_{| \Gr}) \ \subset \ X$, we conclude that $\Pi(\Gr,(\tilde C_X)_{(r)}) \ \subset \ X$ for all $r > 0$,
 so that $\Pi(\bG)_{C^\bu} \ \subset \ X$ and thus $C^\bu$ is an object of $\tilde \cC_X$.
 \end{proof}     

  \vskip .1in

 We now provide an analogue of Theorem \ref{thm:HPS}
for bounded complexes of arbitrary $\bG$-modules.  In this context, the classification is of $(r)$-complete localizing
subcategories of $StMod(\bG)$ replaces the classification of $(r)$-complete tensor ideals of $stmod(\bG)$.
Our proof is heavily dependent upon the classification of localizing subcategories of
 $StMod(\Gr)$ given in \cite{BIKP}.

\vskip .1in

\begin{thm}
\label{thm:HPS-Mod}
  Let $\bG$ be a linear algebraic group.   The correspondences 
$$  X \ \mapsto \ \tilde\cC_X,  \quad \quad  \tilde \cC \ \mapsto \ \Pi(\bG,\tilde\cC)$$
restrict to give mutually inverse bijections
    $$\{ locally \ StMod(\bG){\text-}realizable \ subsets \ X \subset \Pi(\bG) \} $$ 
    
   $$ \longleftrightarrow \quad \quad \{ (r){\text-}complete, \ localizing \ subcategories \  \tilde \cC \subset StMod(\bG) \}.$$
\end{thm}

\begin{proof}
The proof that the evident inclusion $\Pi(\bG,\tilde \cC_X) \subset X$ 
is a bijection if $X$ is a locally  $StMod(\bG)$-realizable subset of $\Pi(\bG)$ is a repetition of the argument 
given in the first part of the proof of Theorem \ref{thm:HPS}.
   
        To complete the proof, we must show that the natural inclusion 
 $\tilde\cC \ \subset \ \tilde \cC_{\Pi(\bG,\tilde\cC)}$ is a bijection for any $(r)$-complete 
 localizing subcategory $\tilde \cC \ \hookrightarrow \ StMod(\bG)$.
Consider some $E^\bu$  which is an object of $\tilde\cC_{\Pi(\bG,{\tilde\cC})} \ \subset \ StMod(\bG)$.    
Since $\tilde \cC$ is $(r)$-complete, it suffices to prove that $E^\bu_{|\Gr}$ is an object of $(\tilde \cC)_{(r)}^\oplus$
for every $r > 0$.  The classification of localizing tensor ideals given in \cite[Thm10.1]{BIKP} tells
us that 
$$(\tilde\cC)_{(r)}^\oplus \quad  =  \quad \tilde \cC_{\Pi(\Gr,{(\tilde\cC)_{(r)}^\oplus})} \ \subset \ StMod(\Gr),$$
so that it suffices to prove that $\Pi(\Gr)_{E^\bu_{|\Gr}} \ \subset \ \Pi(\Gr,(\tilde\cC)_{(r)}^\oplus)$ for 
all $r > 0$.
To prove this inclusion, observe that our assumption on $E^\bu$ implies that
$$\Pi(\Gr)_{E^\bu_{| \Gr}} \ = \ \Pi(\bG)_{E^\bu} \cap \Pi(\Gr) \ \subset \ \Pi(\Gr,\tilde \cC_{|\Gr}) \ 
\ \subset \ \Pi(\Gr,(\tilde\cC)_{(r)}^\oplus)$$
for all $r > 0.$
\end{proof}

\vskip .2in


     \section{Questions and challenges}
     \label{sec:questions-challenges}
     
     Our first remark contrasts our theory of supports $M \mapsto \Pi(\bG)_M$ 
     with various cohomological constructions for $\bG$.
     
     \begin{remark}
    \label{rem:contrast-coh}
    The (rational) cohomology $H^*(\bG,k)$ is invariant under the (conjugation)
    action of $\bG$, whereas the action of $\bG$ on $\Pi(\bG)$ is typically non-trivial.  Thus, it is no surprise that if
    $\bG$ is semi-simple, then $H^*(\bG,k)$ vanishes in positive degrees and $H^*(\bG,M)$ is 
    finite dimensional for any finite dimensional $\bG$-module \cite[Thm2.4]{CPSvdk}.  For $\bG$ unipotent, 
    $H^*(\bG,k)$ does not vanish, but the invariance property implies that  $\Spec H^\bu(\bG,k)$
    is ``too small" to capture much information about $Mod(\bG)$ if $\bG$ is not commutative.
    
    	On the other hand, if $G$ is a finite group scheme, then $|G| \equiv \Spec H^\bu(G,k)$ leads to a
    ``good" cohomological support theory for $mod(G)$, namely $M \ \mapsto \ |G|_M$, 
    where $|G|_M$ is defined as the variety of the annihilator ideal of the $H^\bu(G,k)$-module $Ext_G^*(M,M)$.
    This suggests one might consider for an algebraic group $\bG$ the inverse system $\{ H^\bu(\bG_{(r)},k), r > 0 \}$
    associated to the hyperalgebra $\varinjlim_r k\bG_{(r)}$.
    However, computations in \cite{F5} indicate that
    $\varprojlim_r H^*(\bG_{(r)},M)$ as a $\varprojlim_r H^*(\bG_{(r)},k)$-module is not useful
    in considering $Mod(\bG)$ even in the case of $G$  unipotent.  
    
    Associating a useful cohomology theory for $Mod(\bG)$ using the ``pro-object of finite group schemes"
    $ \cdots \to \bG_{(r+1)} \to \bG_{(r)} \to \cdots $ remains a challenge.
    
   \end{remark}
   
  \vskip .1in
    
    \begin{remark}
    \label{rem:realizable} 
A major challenge is to establish criteria for subsets $X \ \subset \ \Pi(\bG)$ to be realizable
as $\Pi(\bG)_M$ or $\Pi(\bG)_{C^\bu}$.  A better understanding of the support of $\bG$-modules 
obtained by inducing $H$-modules to $\bG$ modules for various closed subgroups of 
$H \ \subset \bG$ would be valuable when investigating this realizability challenge.

A specific challenge is to search for examples of $\bG$-equivariant closed subsets
of $\Pi(\bG)$ which can not be realized as $\Pi(\bG)_{C^\bu}$ for specific groups $\bG$.
\end{remark}
    
    \vskip .1in
    
    \begin{remark}
    \label{rem:balmer}
    It is natural to ask about the Balmer spectrum (see \cite{Bal}) of the triangulated category $stmod(\bG)$
     in light of Theorem \ref{thm:HPS}.  New methods of investigating the Balmer
     spectrum in this non-Noetherian setting are are presumably needed.  In particular, there are typically
    infinitely many isomorphism classes of finite dimensional irreducible $\bG$-modules.
    
    As a start, one might investigate the {\it prime}, $(r)$-complete, tensor ideals of $stmod(\bG)$.  
    These should correspond to ideals of the form $\cC_X \ \subset \ stmod(\bG)$ as
    $X$ varies over subspaces of the form $\Pi(\bG)_{C^\bu}$ which can not be written a non-trivial union of 
    the form $\Pi(\bG)_{D^\bu} \cup \Pi(\bG)_{E^\bu}$ for complexes $C^\bu, \ D^\bu, \ E^\bu$ of
    bounded complexes of finite dimensional $\bG$-modules.   
     \end{remark}
     
    \vskip .1in
      \begin{remark}
    \label{rem:Jordan-function}
    Let $(\bG, \cE)$ be an algebraic group of exponential type.  The support theory 
    $M \ \mapsto \ \bP \cC(\bG)_M \ \simeq \ \Pi(\bG)_M$  extracts 
    minimal information from the 
    ``local operators" $(\cE_{\Lambda_r(\ul B)})_*(u_{r-1})$ 
    and \ $\sum_{s \geq 0}^{r-1} (\cE_{B_s}(u_s)$ at points of $\cC_r(\bG)$.  
    
    For example, we only use the zero locus of the stable Jordan type function of Definition \ref{defn:geom-Cr}. 
 One could ask for restrictions of values of this stable Jordan
 type function for various classes of $\bG$-modules.   As in \cite{FP33},  
    one could  formulate refinements of our support theory $M \ \mapsto \ \bP \cC(\bG)_M$
    using this stable Jordan type function.     
 \end{remark}

\vskip .1in
\begin{remark}
\label{rem:scheme}
 Let $(\bG, \cE)$ be an algebraic group of exponential type.
We do not utilize the scheme structure of $\cC_r(\bG)$, but only the topological
space of scheme-theoretic points.  The scheme structure was used in 
\cite{FP3} and subsequent papers to associate sheaves and vector bundles to 
representations of infinitesimal finite group schemes.  Can this scheme structure be 
used to obtain more information about $\bG$-modules not seen by restrictions
to the family $\{ \Gr \}$ of Frobenius kernels of $\bG$?  
\end{remark}
 
   \vskip .1in

 \begin{remark}
 \label{rem:formal}
 One can consider formal 1-parameter subgroups 
 for a linear algebraic group of exponential type $(\bG,\cE)$.  Namely,  to an
 infinite sequence $\widehat{\ul B} = (B_0,\ldots,B_n,\ldots)$ of pair-wise commuting $K$-points of
 $\cN_p(\fg)$ and a $\bG$-module $M$, one can associate a $p$-nilpotent operator 
 ${\cE_{\widehat{\ul B}}:} M_K \ \to \ M_K$ generalizing the operators we have considered in
 this paper.
 
 At present we lack the technology to use these formal 1-parameter groups to investigate the
 derived category $D^b(Mod(\bG))$.  An appealing feature of such formal 1-parameter groups is 
 that they provide information not seen by $\{ \Gr \}$.
  \end{remark}
    
    \vskip .1in
     
 \begin{remark}
 \label{rem:challenges}  
 In other contexts in which there is an analogous classification theorem, 
knowledge of ``realizable" closed subsets provides insight into the collection of thick tensor ideals 
of a ``representation category."
For $\cC \ \subset \ stmod(\bG)$, could one use some knowledge
of $(r)$-complete, thick tensor ideals to determine interesting classes of finite dimensional $\bG$-modules?
\end{remark}

    \vskip .2in


    \end{document}